\documentclass[a4paper]{article}
    \usepackage[T1]{fontenc}
    \usepackage{mathtools}
    \usepackage{amssymb}
    \usepackage{amsthm}
    \usepackage{authblk}
    \usepackage{algorithm}
    \usepackage{algpseudocodex}
    \usepackage{mathtools}
    \usepackage{enumerate}
    \usepackage{bm}
    \usepackage{ulem}
    \usepackage{bbding}
    \usepackage{xfrac}
    \usepackage{subcaption}
    \usepackage{pgfplots}
    \usepackage{titlesec}
    \usepackage{rotating}
    \usepackage{comment}
    \usepackage{cite}
    \usepackage[final]{showlabels}
    \usepackage[textsize=footnotesize,textwidth=2cm]{todonotes}
    \usepackage[hidelinks,hypertexnames=false]{hyperref}
    \usepackage[capitalise]{cleveref}
    \usepackage{xcolor}
    \usepackage{xargs}
    \usepackage{colortbl}
    \usepackage{multicol}

    \newcommand{\defeq}{\vcentcolon=}
    \newcommand{\defeqr}{=\vcentcolon}
    \newcommand*{\tran}{^{\mkern-1.5mu\mathsf{T}}}
    \newcommandx{\axtd}[2][1=]{\todo[linecolor=orange,backgroundcolor=orange!40,bordercolor=orange,#1]{#2}}
    \newcommandx{\dmtd}[2][1=]{\todo[linecolor=magenta,backgroundcolor=magenta!40,bordercolor=magenta,#1]{#2}}

    \DeclareMathOperator{\erf}{erf}
    \DeclareMathOperator{\erfc}{erfc}

    \title{A flexible and robust approach to univariate Gaussian splitting using parameterised
    Gaussian mixtures}
    \author{Dmitry Mikhin and Athena Xiourouppa}
    \affil{Acacia Systems Pty Ltd}
    \date{\today}

    \pgfplotsset{compat=1.18}

\titleformat{\paragraph}
{\normalfont\normalsize\bfseries}{\theparagraph}{1em}{}
\titlespacing*{\paragraph}
{0pt}{3.25ex plus 1ex minus .2ex}{1.5ex plus .2ex}

\titleformat{\subparagraph}
{\normalfont\normalsize\bfseries}{\thesubparagraph}{1em}{}
\titlespacing*{\subparagraph}
{0pt}{3.25ex plus 1.5ex minus .2ex}{2.5ex plus .2ex}

\makeatletter
\@addtoreset{section}{part}
\makeatother

\titleformat{\subparagraph}
{\normalfont\normalsize\bfseries}{\thesubparagraph}{1em}{}
\titlespacing*{\subparagraph}
{0pt}{3.25ex plus 1ex minus .2ex}{1.5ex plus .2ex}

\makeatletter
\newcommand{\raisemath}[1]{\mathpalette{\raisem@th{#1}}}
\newcommand{\raisem@th}[3]{\raisebox{#1}{$#2#3$}}
\makeatother

\begin{document}

\maketitle

\pgfmathdeclarefunction{gauss}{3}{%
  \pgfmathparse{#3/(#2*sqrt(2*pi))*exp(-((x-#1)^2)/(2*#2^2))}%
}

\begin{abstract}
    We consider approximation of a Gaussian distribution with a mixture of homoscedastic Gaussians of
    smaller variance. The solution is obtained by minimising the $L^2$ norm between the original
    Gaussian and the mixture, which is parameterised to reduce the complexity of the
    optimisation problem.
    The developed technique is straightforward, sufficiently robust and yields
    Gaussian Mixtures that rapidly approach the original function as the number of mixands is
    increased. The proposed solution is examined for multiple special cases of input parameters
    resulting in further simplifications. Extension of the proposed method for approximating
    non-Gaussian distributions is discussed.
\end{abstract}

\section{Introduction}

Application of finite mixture models for analysis of complex data is ubiquitous in statistical
literature (\textit{e.g.} \cite{McLachlan:2019:Finite,Richardson:1997:Bayesian} and numerous references within). If
component densities are taken from a parametric family, such as the multivariate
normal, mixture models combine the benefits of parametric and nonparametric
approaches to statistical estimation.
The expectation-maximization (EM) algorithm
\cite{Dempster:1977:ML} provides a practical way of fitting mixture models to experimental data. In tracking applications, however, finite mixtures
of normal components may arise in a different context.

The classic Kalman Filter (KF) \cite{Kalman:1960:New} and its extensions \cite{Maybeck:1979:Stochastic,Sarkka:2013:Bayesian} describe the state of the system as
a Gaussian distribution. This state is propagated forward in time by iteratively repeating two steps:
time update (prediction) and measurement update. If the
time update equations are non-linear, the exact
predicted state would not be a Gaussian distribution. Then
the expected measurement is derived from the predicted state using measurement equations that also could be
non-linear. Thus, neither the predicted state nor the predicted measurement distributions are
Gaussian. At best, they
can be \textit{approximated} by Gaussians. The Unscented KF (UKF) \cite{Wan:2001:Unscented} and
Cubature KF (CKF) \cite{Arasaratnam:2009:Cubature} are examples of practical algorithms that provide
good approximations, but
the more non-linear the prediction or measurement equations are, the stronger the deviation from normality could be.

The strength of non-linearity has to be quantified with respect to the state distribution itself.
The prediction functions are usually continuously differentiable, and therefore,
can be linearly approximated within some intervals of their variables. If the
input distributions are sufficiently narrow, they may fit within the linearity limits, the
updated distributions would be close to normal, and UKF / CKF would describe them well.
Problems occur when the equations are substantially non-linear over the spread of the
original distributions.

One possible solution is to ``split'' the original Gaussian distribution
and approximate it as a sum of
several Gaussians with smaller covariances \cite{Alspach:1972:Nonlinear,DeMars:2013:Entropy,Faubel:2009:Split,Vittaldev:2016:Multidirectional}. The narrower distributions would be less
affected by non-linearity and their approximations by UKF / CKF algorithms would be more accurate.
In this approach, the state and measurement distributions are represented as Gaussian Mixture models, and the
resulting tracking method is known as the Gaussian Sum Filter (GSF) \cite{Alspach:1972:Nonlinear,
Leong:2013:Gaussian}.

GSF needs to determine when a
split is appropriate, find the optimal direction(s) of splitting,
and then apply the actual splitting approximation \cite{Zanetti:2025:Uncertainty,Sun:2023:Hybrid,Chen:2024:Enhanced}. Multiple approaches for the first two
subtasks are reviewed in \cite{Zanetti:2018:Novel,Kulik:2025:Nonlinearity}. This
paper is focused on the last step: we present an algorithm for splitting a
univariate Gaussian distribution into homoscedastic Gaussian Mixtures; such univariate methods serve
as the
foundation for multivariate splitting algorithms
\cite{DeMars:2013:Entropy,Vittaldev:2016:Space,Sun:2023:Hybrid}.
Our approach is based on minimising the $L^2$ norm between the original
distribution and the mixture model. We carry out asymptotic and numeric analysis of the
obtained approximation, and determine the parameters for yielding the $L^2$ mismatch within a given
threshold.

The univariate splitting problem is summarised in the opening \cref{sec:ReqAndCons}.
We review the existing splitting methods in \cref{sec:LiteratureReview}, and identify some issues
that motivate our approach. The new algorithm and its asymptotic analyses for various special cases
are given in \cref{sec:Split1D}, followed by the computational results in \cref{sec:Numerics:1D}.
The concluding \cref{sec:Discussion} summarises the results and discusses
possible extensions.

\section{Requirements and constraints}\label{sec:ReqAndCons}

The univariate problem is deceptively simple: given the standard Gaussian distribution
$\widetilde{\mathcal{N}} (x; 0, 1)$\footnote{Without any loss of generality, we aim to split the
standard Gaussian; all other univariate cases are covered by a linear change in coordinate,
\textit{i.e.,} offset and scaling.}, we aim to find the means $\mu_m$, variances $\sigma_m^2$, and
weights $w_m$ that provide the best approximation of the form
\begin{equation}
    \label{eq:Approx}
    \widetilde{\mathcal{N}}(x; 0, 1) \approx \widetilde{\mathcal{Q}}(x) = \sum_{m = 0}^{M - 1} w_m \mathcal{N}(x; \mu_m, \sigma_m^2) \text{.}
\end{equation}
Here $\mathcal{N}(x; \mu, \sigma^2)$ is the PDF of a univariate normal distribution with the
mean $\mu$ and variance $\sigma^2$.
The ``best'' approximation could be the one that minimises the
Kullback-Leibler (KL) divergence \cite{Cover:2006:Elements}, or the Hellinger distance
\cite{Gibbs:2002:Choosing}, or some other measure of difference between $\widetilde{\mathcal{N}}$ and
$\widetilde{\mathcal{Q}}$; $\mu_m$, $w_m$, and $\sigma^2_m$ are found as the solutions of
the obtained optimisation problem. Alternatively, we may require that the
approximation has the same first, second, and some higher-order moments as the
original distribution; these moment-matching equations are solved for the parameters of $\widetilde{\mathcal{Q}}$.

Regardless of the selected method, the approximation always has a trivial solution:
set $w_0 = 1$, $\mu_0 = 0$, $\sigma_0 = 1$, $w_m = 0$, $\forall m \ne 0$
and \cref{eq:Approx} becomes an equality; the best Gaussian Sum (GS) approximation
to a Gaussian function is that function itself, without any splitting. We must
restrict the search space to exclude this undesired ``absolute best'' result. Also
the solution has to satisfy a few additional requirements. The
weights must be normalised as
\begin{equation}
    \label{eq:WNorm}
    \sum_{m = 0}^{M - 1} w_m = 1 \text{,}
\end{equation}
because the right-hand side of \cref{eq:Approx} must be a PDF of a probability distribution.
Second, we generally demand that all weights are positive, as
the PDF may
not take negative values.\footnote{This condition is sufficient but not strictly required; only the sum of
all the terms in \cref{eq:Approx} must be non-negative, which theoretically can be satisfied even if some weights
are negative.}

\section{Literature review}\label{sec:LiteratureReview}

\subsection{Existing splitting methods}\label{sec:Existing}

A number of splitting algorithms for GS filters were developed
starting with \cite{Sorenson:1971:Recursive}, who
approximated uniform and gamma, but not Gaussian, distributions. Regrettably, splitting can be
viewed as a technicality and covered briefly, especially for minor variations of the same
techniques.
Our motivation for splitting is to make the Gaussian components narrow enough so that
linearisation works. Thus, the following literature review summarises the proposed methods
and the achieved variance reduction, if reported.

The twin problem of Gaussian splitting is Gaussian reduction or merging, which aims to approximate a
large GS with a smaller one. We touch upon some reduction algorithms that can be
re-purposed for splitting.

\subsubsection{Moment-matching and related methods}

Simple, closed-form methods for splitting into two or three components are proposed in
\cite{Faubel:2010:Further}. Both algorithms employ a moment-matching approach: the mean
and the covariance are matched exactly; the two-term solution approximates the fourth moment, while the three-term
solution matches the fourth moment exactly and approximates the sixth.
Instead of solving this moment-matching problem directly, the authors derive
parameterised expression for the means, covariances, and weights, and then select the
parameter values that
``give a good trade-off between displacement of components and accuracy of approximation'' \cite[p.~851]{Faubel:2010:Further}.
For the two-term approximation, both mixands have the standard
deviation of $\sigma \approx 0.866$ for the recommended parameter value.
In the three-term case, the standard deviation is even larger, at $\sigma
\approx 0.957$.
The two-term variant is later used in
\cite{Huber:2011:Adaptive}.

Building on the results of \cite{Faubel:2010:Further}, \cite{Leong:2013:Gaussian} apply the two-term
algorithm recursively to obtain a $2 \times 2 = 4$-term split.
Double application of the two-term
splitting reduces the standard deviation to $\sigma
= 0.75$.
However, the resulting GS
has the same weights for all terms, which compromises the approximation of the
original distribution shape. If pair-wise splitting of existing components is applied,
all other components should be adjusted to minimise some global measure of approximation error.

A multivariate moment-matching splitting algorithm with an arbitrarily large (odd) number of mixands
is proposed in \cite{Leutnant:2011:Versatile}. Mixand are homoscedastic with a constant
separation vector between the means. The method allows simultaneous splitting along
multiple directions. General conditions are established to ensure non-negative definiteness of the
mixands. No closed-form solution is found for the weights and they are assigned
\textit{a priori}.
Even matching fourth-order moments is problematic, as it conflicts with the requirements of
positive weights and positive-definite mixands covariances. The univariate
approximation is not considered, so we derive it in \cref{sec:Leutnant:1D} and find the maximal
possible reduction of the standard deviation to be $\sigma \approx 0.272$, achieved in the limit
of $M \to \infty$.
\cite[Eqs. 33-34]{Jiang:2020:Gaussian} use a closed-form
three-component moment-matching splitting attributed to
\cite{Leutnant:2011:Versatile}.

The splitting algorithm by \cite{Raitoharju:2015:Binomial} is based on
the observation that the sequence of standardised binomial distributions converges weakly to a
standard normal distribution as the number of trials increases. Then the GS
approximation is constructed using component means and weights from scaled binomial distributions. The number of
components can be arbitrarily large; examples in \cite{Raitoharju:2015:Binomial} use splittings with up to $M = 512$ terms.
The approximation is evaluated only in the transformed space, after passing
through the non-linear function; no direct comparison to the original distribution is given.

A technique for splitting a multivariate Gaussian into two by scaling covariance matrix eigenvectors
(alternatively, columns of its Cholesky decomposition) is proposed in \cite{Zhang:2003:EM} and reused for tracking applications by \cite{Terejanu:2011:SplitMerge}. The
solution is expressed in closed form and preserves the first and second moments. The variance
reduction is not discussed, but cannot be significant for the one-to-two split. Other mixture
components are not adjusted as part of the algorithm, so, the approximation of the overall distribution shape
may be compromised.

\subsubsection{Multivariate optimisation methods}

A general method for approximating density functions by solving
a system of differential equations for the approximation parameters is described in \cite{Hanebeck:2003:Progressive}. However, the
presented technique for splitting a Gaussian function
differs from the general approach. First, the distribution is split into two
components with equal variances, then those variances are reduced while adapting the component
weights and means (presumably to minimise the deviation from the true function; no details are
given). When the standard deviations ``cannot be reduced any further'' \cite[p.~8]{Hanebeck:2003:Progressive},
the splitting is repeated, yielding four (then eight, sixteen, and so on)
homoscedastic components, and the variance reduction is repeated ``until the desired number
\ldots of mixture components is achieved'' \cite[p.~8]{Hanebeck:2003:Progressive}.
In addition to some complexity inherent in
recursive splitting, when the same approach is used in \cite{Schrempf:2005:Optimal} for simplifying
products of Gaussian mixtures, its results are inferior to the simpler method of
\cite{Schieferdecker:2009:Gaussian}.
The only example presented in \cite{Hanebeck:2003:Progressive} has four components
(figs.~3 and 4) and achieves variance reduction of
$\sigma \approx 0.518$. Later
\cite{Huber:2008:Entropy} use a four-component approximation with a reference to
\cite{Hanebeck:2003:Progressive};
the same splitting is applied by
\cite{Straka:2016:Directional}, without details.

The method of \cite{DeMars:2013:Entropy} uses the beta divergence \cite{Cichocki:2010:Families}
as the difference measure between the original Gaussian and its approximation; the selected beta
divergence parameters, however, simplify it to the equivalent of the squared $L^2$
distance. The approximation is homoscedastic. The solution is
obtained by multivariate optimisation with regularisation: the function to be minimised is the
$L^2$ norm squared plus the weighted mixand variance. This
automatically avoids the undesired ``best possible'' solution of a single Gaussian. On the downside,
the user-specified weight of the $\sigma^2$ term provides only indirect control of the optimal
variance. The parameters of the three- and five-component splittings are given in \cite[tables~1
and~2]{DeMars:2013:Entropy}.
The three-component variant has $\sigma \approx 0.672$, reduced to $\sigma \approx 0.442$ for
five components.
In \cite{LeGrand:2022:Split}, this method is used in a multivariate
context for recursive splitting along multiple directions.

A more restricted version of the \cite{DeMars:2013:Entropy} approach is used in \cite{Kulik:2025:Nonlinearity}. The mixand
variance is additionally constrained to preserve the variance of the original distribution, while
the means are uniformly spaced, leaving only the weights as the free parameters.
More constraints are imposed to ensure that the variance remains positive.
The presented multivariate examples apply the univariate algorithm
recursively along several directions, resulting in a total of nine (in example A, recursion depth of
two), 81 (example B, depth of four), and 27 (example C, depth of three) mixands, suggesting that the
univariate splitting in all these cases uses only three mixands.

The optimisation approach of \cite{Tuggle:2020:Dissertation,Zanetti:2025:Uncertainty} uses the KL divergence between the
standard Gaussian and its GS approximation as the cost function for minimisation. The components are
homoscedastic, but not equidistant: their means are set deterministically as the
centres of equi-probable regions of the distribution. The variance is then fixed by second moment
matching, leaving only the weights as the optimisation arguments. Symmetry considerations and
the normalisation condition by \cref{eq:WNorm} reduce the number of unknowns to just one
for three- and four-term mixtures, and to two for five- and six-term mixtures. The splittings are
illustrated in \cite[figs.~B.2-B.5]{Tuggle:2020:Dissertation} and \cite[fig.~1]{Zanetti:2025:Uncertainty}, but their parameters are given
explicitly only for the three-term case \cite[table~4.3]{Tuggle:2020:Dissertation}, which yields the
variance reduction of $\sigma \approx 0.784$. This is larger than for the three-term
approximation of \cite{DeMars:2013:Entropy}, and it may explain the preference for the latter in the
examples of \cite{Tuggle:2020:Dissertation}.
The splittings from \cite{DeMars:2013:Entropy} are also
used earlier in \cite{Tuggle:2018:Automated}.
The choice of KL divergence as the cost function is attractive from the information-theoretical point of
view, but presents substantial implementation difficulties because this measure is not available in
closed form for comparing Gaussian and GS distributions.

In a major extension of univariate splitting capabilities, \cite{Vittaldev:2016:Multidirectional}
apply advanced multidimensional optimisation techniques to obtain the best approximation in terms
of $L^2$ distance for an odd number of mixands up to 39. The $L^2$ measure is used because both the
norm and its derivatives can be found in closed form. Still, this cost function is characterised as
``highly nonlinear with many local minima'' \cite[p.~92]{Vittaldev:2016:Multidirectional}. The upper limit of 39 terms is attributed to precision
limitations, despite using quadruple-precision floating-point arithmetic. The
variance of homoscedastic mixands is given by pre-defined rules as a low negative power of $M$, the
number of terms in the mixture. Three possible ``rules'' for $\sigma^2(M)$ are considered, and the
presented $L^1$ and $L^2$ norms of the approximation error radically depend on the selected rule.
When these univariate splittings are applied in the multivariate context, the results
are sensitive to the variance rule, and the optimal rule apparently depends on the multivariate
problem at hand. The developed splittings are published online; one seven-term variant is detailed
in the paper. The method is later used by \cite{Vittaldev:2016:Space,Sun:2023:Hybrid,Chen:2024:Enhanced}.

The standard deviation given by the \cite{Vittaldev:2016:Multidirectional} library of splittings
depends on the number of terms and the selected $\sigma^2(M)$ rule. For the highest number of terms
possible for each rule it is $\sigma \approx 0.160$ in the $1 / M$ rule for $M = 39$, $\sigma
\approx 0.299$ in the $(1 / M)^{\sfrac{3}{4}}$ rule for $M = 25$, and $\sigma \approx 0.508$ in the
$(1 / M)^{\sfrac{1}{2}}$ rule for $M = 15$. The significant reduction for the inverse-$M$ rule is
questionable, as the mixands are so narrow that they can nominally cover only about $\sigma \times M
= 1$ standard deviation of the original distribution. Out of the three rules, this
one has the worst $L^1$ and $L^2$ approximation errors.

Considering now the publications on Gaussian merging, which aim to reduce the number of terms in
a GS of many components, while maintaining the overall shape of the distribution.
Some of their observations
equally apply to the splitting task.
Upon reviewing several
measures of difference between the original and reduced forms, \cite{Williams:2003:Cost} settle on
the Integrated Square Difference, which is the same as the $L^2$ norm squared. The actual
reduction involves solving a non-linear optimisation problem, and false local minima are an
issue, such that selecting a good initial point for optimisation is critical. The techniques of
\cite{Schieferdecker:2009:Gaussian, Chang:2010:Scalable}, evaluated in depth in
\cite{Crouse:2011:Gaussian}, use the same metric (under the name of Integral Squared Error)
for mixture-reduction purposes. These papers also report the problem of multiple local
minima, so, by corollary, we may expect it for splitting algorithms that use multivariate
optimisation.

\subsubsection{Sigma-point methods}

In \cite{Vishwajeet:2018:Adaptive}, splitting is integrated directly into the multivariate UKF state update equations.
The components are homoscedastic and have equal weights.
Under such assumptions, the component parameters are found from
conservation of weight, mean, and the second moment by the GS approximation. No direct comparison
with the original distribution is presented.

Similarly, \cite{LeGrand:2023:Cislunar} propose a splitting method based on Cholesky
square root factors, which is intricately weaved into the process of propagating a multivariate
Gaussian component through a non-linear transformation by sigma-point sampling.
No direct assessment of the splitting accuracy is provided.

\subsection{A recap}\label{sec:Recap}

In summary, the three most common approaches to splitting are moment-matching,
exploiting the inner workings of the sigma-based Kalman filters, and
multivariate minimisation of the approximation error.

Moment-matching techniques are usually restricted to low moment orders.
Conservation of higher-order moments results in correspondingly high-order
equations for the mixand parameters, which quickly become unsolvable in closed form, and
likely hard to analyse numerically. Some additional \textit{a
priori} assumptions are used to obtain the solution. Non-negativity constraints on
weights and variances might preclude exact moment matching.

The higher-order moments are instrumental in describing the distribution behaviour at
the tails. In tracking problems, however, the tails are of less interest, as
the target is not expected to be at infinity, even with an infinitesimally small probability.
Thus, the higher-order moments are of less concern
and shape matching at the core of the distribution may be preferable.

The sigma-point methods are a special variant of moment-matching, integrated with filter updates.
Only lower-order moments are considered and no analysis is available for the approximation accuracy
with respect to the original distribution.

The optimisation approach usually results in constrained problems. The $L^2$
norm is the most common cost function. The problem of multiple local minima is often reported.
Many of optimisation approaches
use at least some elements of moment matching, such as conservation of the first and second moments.

The majority of the reviewed methods present splittings for low to intermediate numbers of terms.
The reported reduction in standard deviation is moderate for nearly all methods. The only exception
is the technique of \cite{Vittaldev:2016:Multidirectional} with the inverse-$M$ rule, but it poorly
approximates the original Gaussian. Therefore, using the presented
splitting libraries, we can expect only modest suppression of the non-linearity effects across the
spread of the distribution.

A common theme \cite{Faubel:2010:Further, Jiang:2020:Gaussian,
Kulik:2025:Nonlinearity, Leutnant:2011:Versatile, Tuggle:2020:Dissertation} is to start with the
offset between the adjacent component means as the primary parameter, and then use second moment
matching to express the homoscedastic mixand variance through the weights and the offset. Then
other conditions are applied to find the weights, so the variance can be computed, but it
might end up insignificantly reduced compared to the original distribution or even negative.

This approach to the problem seems to go the wrong way. Our primary values of interest are the
component variance $\sigma^2$, because it determines the impact of non-linearity, and the number of
splitting terms $M$, because it determines the computational load and
complexity of the algorithm. It is arguably better to specify the design constraints in terms of
$\sigma$ and $M$ and search for the optimal component separation and weights.

\section{The new algorithm}\label{sec:Split1D}

We seek the best approximation in \cref{eq:Approx}, subject to constraints, such that the
undesired solution of a single non-zero component is excluded. As argued above, the
superior way for describing the desired approximation starts with specifying the design constraints
on the component variance and the number of terms. As an immediate bonus, limiting or directly
specifying the component variance already eliminates the undesired solution. To approximate
the distribution shape, we choose the most direct metric, the squared $L^2$ norm of the
mismatch, defined as
\begin{equation}
    \label{eq:L2}
    L^2 = \int_{-\infty}^{+\infty}
        {\left[ \widetilde{\mathcal{N}}(x; 0, 1) - \widetilde{\mathcal{Q}}(x) \right]}^2 \! dx \text{.}
\end{equation}
We then simplify the problem by using homoscedastic components with variances
$\sigma_m^2 = \sigma^2$ and equidistant means $\mu_m = m h + h_M$, where $h$ is the step
between the means, and $h_M$ is the initial offset. From the desired symmetry of the approximation
we get $h_M = 0$ when $M$ is odd and $h_M = h / 2$ when $M$ is even.
The remaining approximation parameters are $h$ and $w_m$, to be found by minimising $L^2$.
The symmetry allows to reduce the number of these parameters even further, but
the specifics depend on $M$ being odd or even, and are considered in the respective
subsections. For now, we make another general observation. The conditions for the
$L^2$ minimum through the partial derivatives are
\begin{equation}
    \label{eq:Partial}
    \frac{\partial L^2}{\partial w_m} = 0 \text{, for } m = 0, \ldots, M - 1 \text{, and }
    \frac{\partial L^2}{\partial h} = 0 \text{.}
\end{equation}
It follows from \cref{eq:L2} that $L^2$ is a quadratic form in $w_m$, and therefore, the first $M$
of equations~\labelcref{eq:Partial} form a linear system for the weights $w_m$.
Solving this system allows expressing the weights as functions of $h$. Then $L^2$ becomes a
univariate function of $h$, and the last of \cref{eq:Partial} becomes a univariate non-linear
equation for $h$. To solve it, many powerful
numerical algorithms are immediately available.

\subsection{Odd case}
\label{sec:1D:Odd}

We begin with the detailed solution for the case when the total number of terms is odd.
Changing the
notation slightly, $M \Rightarrow 2 M + 1$, the component index $m$ in
\cref{eq:Approx,eq:WNorm,eq:L2} would run from $-M$ to $M$. By the symmetry of approximation we then have
$w_{-m} \equiv w_m$, $\forall m = 1, \ldots, M$. The normalisation condition from \cref{eq:WNorm}
gives
\begin{equation}
    \label{eq:w0:Odd}
    w_0 = 1 - 2 \sum_{m = 1}^{M} w_m \text{.}
\end{equation}
Thus, by symmetry and normalisation the number of unknowns is reduced by half. Now, substituting
\cref{eq:w0:Odd} into \cref{eq:L2}, we obtain
\begin{equation}
    \label{eq:L2:Odd}
    L^2 = \int_{-\infty}^{+\infty}
        {\left[ \widetilde{\mathcal{N}} - \mathcal{N}_0 -
                \sum_{m = 1}^{M} w_m \bigl( \mathcal{N}_m - 2 \mathcal{N}_0 + \mathcal{N}_{-m} \bigr) \right]}^2 \! dx \text{.}
\end{equation}
We omit the arguments to the various normal PDFs, and use a
shorthand notation $\mathcal{N}_\alpha \defeq \mathcal{N}(x; \alpha h, \sigma^2)$. To
utilise the quadratic dependency of $L^2$ on $w_m$, we introduce basis functions
\begin{equation}
    \label{eq:Base:Odd}
    f_m \defeq \mathcal{N}_m - 2 \mathcal{N}_0 + \mathcal{N}_{-m} \text{,}
\end{equation}
and also denote
\begin{equation}
    \label{eq:RHS:Odd}
    F \defeq \widetilde{\mathcal{N}} - \mathcal{N}_0 \text{.}
\end{equation}
Although the $f_m$ are not orthogonal, they are not linearly dependent either, except for the
degenerate case of $h = 0$. Using $f_m$ and $F$, \cref{eq:L2:Odd} simplifies to
\begin{equation}
    \label{eq:L2:Quad}
    L^2 = \int_{-\infty}^{+\infty} {\Bigl[ F - \sum_m w_m f_m \Bigr]}^2 dx
        = {\lVert F \rVert}_2 - 2 \bm{w}\tran \bm{b} + \bm{w}\tran \mathbf{A} \bm{w}
        = {\lVert F \rVert}_2 - \bm{w}\tran \bm{b} \text{,}
\end{equation}
where $\bm{b} = {[b_1, \ldots, b_M]}\tran$, $\mathbf{A} = \lVert a_{m,k} \rVert$, and
\begin{equation}
    \notag
    \begin{aligned}
        {\lVert F \rVert}_2 &= \langle F, F \rangle = \int_{-\infty}^{+\infty} F^2 dx \text{,} \\
        b_m &= \langle F, f_m \rangle = \int_{-\infty}^{+\infty} F f_m dx \text{,} \\
        a_{m,k} &= \langle f_m, f_k \rangle = \int_{-\infty}^{+\infty} f_m f_k dx \text{.} \\
    \end{aligned}
\end{equation}
The functions $f_m$ and $F$ are linear combinations of Gaussian PDFs, and the product of Gaussian PDFs
is a scaled Gaussian PDF. Thus, using \cref{eq:Base:Odd,eq:RHS:Odd}, we obtain for $b_m$, $a_{m,k}$,
and ${\lVert F \rVert}_2$
\begin{equation}
    \label{eq:cd:Odd}
    \begin{aligned}
        {\lVert F \rVert}_2 &= {(2 \sqrt{\pi})}^{-1} - 2 d_0 + c_{0,0} \text{,} \\
        b_m &= 2 (c_{0,0} - c_{0,m} - d_0 + d_m) \text{,} \\
        a_{m,k} &= 2 (c_{k,m} + c_{-k,m} - 2 c_{0,m} - 2 c_{0,k} + 2 c_{0,0}) \text{,} \\
    \end{aligned}
\end{equation}
where
\begin{equation}
    \label{eq:Gprods}
    \begin{aligned}
        d_\alpha &\defeq \langle \widetilde{\mathcal{N}}, \mathcal{N}_\alpha \rangle
            = \frac{1}{\sqrt{2 \pi (1 + \sigma^2)}} \exp{\Bigl( -\frac{\alpha^2 h^2}{2 (1 + \sigma^2)} \Bigr)} \text{,} \\
        c_{\alpha,\beta} &\defeq \langle \mathcal{N}_\alpha, \mathcal{N}_\beta \rangle
            = \frac{1}{2 \sigma \sqrt{\pi}} \exp{\Bigl( -\frac{{(\alpha - \beta)}^2 h^2}{4 \sigma^2} \Bigr)} \text{.} \\
    \end{aligned}
\end{equation}
Although here the indices $\alpha$, $\beta$ of $c_{\alpha,\beta}$ are integers,
\cref{eq:Gprods} are valid in a general case when $\alpha,\beta \in \mathbb{R}$.

Taking partial derivatives of \cref{eq:L2:Quad} over $w_m$ and solving the obtained linear system as
$\bm{w} = \mathbf{A}^{-1} \bm{b}$, we express $L^2$ as a univariate function of $h$
\begin{equation}
    \label{eq:L2:special}
    L^2 = \mathcal{M}(h; \sigma, M) \text{.}
\end{equation}
It is possible that $\mathcal{M}$ is a known special function, but so far establishing such an
association eluded us. Thus, we could not rely on an existing theory of $\mathcal{M}$ properties
and have to establish them from scratch.
The minimum of $L^2(h)$ is found numerically. Furthermore, we cannot guarantee that
the obtained weights are non-negative; we can only check their signs \textit{post facto}.
Some properties of the weights can be established analytically, in particular, we show that
for certain $h$ the weights are positive or cannot be positive.

\subsubsection{Asymptotic solution for large \texorpdfstring{$h$}{h}}
\label{sec:LargeH:Odd}

In the limit of large steps, $h \to \infty$, the $c_{\alpha,\beta}$ factors tend to zero unless
$\alpha = \beta$, and the $d_\alpha$ factors tend to zero unless $\alpha = 0$. Therefore,
\cref{eq:cd:Odd} simplify to:
\begin{equation}
    \notag
    \begin{aligned}
        b_m &\to 2 (c_{0,0} - d_0)
            &&= \frac{1}{\sigma \sqrt{\pi}} \Bigl( 1 - \frac{2 \sigma}{\sqrt{2 (1 + \sigma^2)}} \Bigr)
            \defeqr b_\infty \text{,} \\
        a_{m,k} &\to 2 (c_{k,m} \delta_{k,m} + 2 c_{0,0})
            &&= \frac{1}{\sigma \sqrt{\pi}} (2 + \delta_{k,m}) \text{,} \\
    \end{aligned}
\end{equation}
for all $m$ and $k$; here $\delta_{k,m}$ is the Kronecker delta symbol. Thus, we have
\begin{equation}
    \notag
    \mathbf{A} \to \frac{1}{\sigma \sqrt{\pi}} (2 \mathbf{J} + \mathbf{I}) \quad \text{and} \quad
    \bm{b} \to b_\infty \bm{j} \text{,}
\end{equation}
where $\mathbf{J}$ is a matrix of all ones, $\mathbf{I}$ is the identity matrix, and $\bm{j}$ is a
vector of all ones. From the normalisation condition in \cref{eq:w0:Odd}, we have $\mathbf{J}
\bm{w} = (1 - w_0) \bm{j} / 2$, and therefore,
\begin{equation}
    \notag
    \mathbf{A} \bm{w} \approx \frac{1}{\sigma \sqrt{\pi}} (2 \mathbf{J} + \mathbf{I}) \bm{w}
                      = \frac{1}{\sigma \sqrt{\pi}} \left( (1 - w_0) \bm{j} + \bm{w} \right)
                      \approx b_\infty \bm{j} \text{.}
\end{equation}
Solving this equation, the large-step approximation for the weights is
\begin{equation}
    \notag
    \bm{w} \approx \Bigl( w_0 - \frac{2 \sigma}{\sqrt{2 (1 + \sigma^2)}} \Bigr) \bm{j} \text{.}
\end{equation}
Substituting it back into the normalisation condition and solving for $w_0$, we find
\begin{equation}
    \label{eq:W:Odd:Inf}
    \begin{aligned}
        \lim\limits_{h \to \infty} w_0
            &= \frac{1}{1 + 2 M} \Bigl( 1 + \frac{4 M \sigma}{\sqrt{2 (1 + \sigma^2)}} \Bigr) \text{,} \\
        \lim\limits_{h \to \infty} w_m
            &= \frac{1}{1 + 2 M} \Bigl( 1 - \frac{2 \sigma}{\sqrt{2 (1 + \sigma^2)}} \Bigr) \text{, } m \ne 0 \text{.} \\
    \end{aligned}
\end{equation}
All asymptotic weights are non-negative, as $\sigma \le 1$. For the $L^2$ norm at large $h$ we
obtain
\begin{equation}
    \label{eq:AsInf:Odd}
    L^2_\infty = \frac{1}{2 \sqrt{\pi}}
        - \frac{2}{\sqrt{2 \pi (1 + \sigma^2)}}
        + \frac{1}{2 \sigma \sqrt{\pi}}
        - \frac{M \sigma \sqrt{\pi}}{1 + 2 M} b^2_\infty
        \text{.}
\end{equation}

\subsubsection{Asymptotic solution for small \texorpdfstring{$h$}{h}}
\label{sec:SmallH:Odd}

Now we explore the opposite limit of $h \to 0$. The details of this long derivation are given
in \cite{Mikhin:2026:Asymptotic} and only a summary is presented here.

We represent $L^2$ at small $h$ as a power series over $h_* = h / \sigma$:
\begin{equation}
    \label{eq:decL2:h0:Odd}
    L^2 = \sum_{j=0}^{\infty} L^2_{2j} h_*^{2j} \text{.}
\end{equation}
The odd terms are zero by symmetry. To find $L^2_{2j}$ we need the corresponding power series for
$\bm{w}$ and $\bm{b}$, and to obtain the former, we need the power series of $\bm{A}$.
Differentiating the elements $b_m$ from \cref{eq:cd:Odd}
with respect to $h_*$ and evaluating at $h_* = 0$, we obtain
\begin{equation}
    \label{eq:decb:h0:Odd}
    \bm{b} = \sum_{k=1}^{\infty} \beta_{2k} \bm{p}_{2k} h_*^{2k} \text{,}
\end{equation}
where
\begin{equation}
    \label{eq:beta:h0:Odd}
    \beta_{2k} =
        \frac{(-1)^k}{k! \, 2^{2k} \sigma \sqrt{\pi}} \left( \frac{2^{k + 1/2} \sigma^{2k + 1}}{(1 + \sigma^2)^{k + 1/2}} - 1 \right) \\
    \text{,}
\end{equation}
and
\begin{equation}
    \label{eq:p2k:Odd}
    \bm{p}_{2k} = \bigl[ 1, 2^{2k}, \dots, M^{2k} \bigr] \text{, } \forall k \ge 1 \text{.}
\end{equation}
The vectors $\bm{p}_{2k}$ for $k \le M$ are independent, although not mutually orthogonal.

Similarly, differentiation of $a_{m, k}$ from \cref{eq:cd:Odd} with respect to $h_*$ at
$h_* = 0$ yields the power series for the matrix $\bm{A}$
\begin{equation}
    \label{eq:deca:h0:Odd}
        \bm{A} = \sum_{n=2}^{\infty} \bm{A}_{2n} h_*^{2n} \text{,}
\end{equation}
where the matrices $\bm{A}_{2n}$ are expressed as
\begin{equation}
    \label{eq:A:Odd}
    \begin{aligned}
        \bm{A}_{2n} = \sum_{j = 1}^{n - 1} \alpha^{(2n)}_{2j,2(n - j)}
                          \left( \bm{p}_{2 j} \otimes \bm{p}_{2(n - j)} \right) \text{,}
    \end{aligned}
\end{equation}
the $\otimes$ symbol denotes the outer product of vectors, and
\begin{equation}
    \label{eq:alpha:h0:Odd}
        \alpha^{(2n)}_{2j, 2(n - j)} =
            \frac{(-1)^n}{n! \, 2^{2n - 1} \sigma \sqrt{\pi}} \binom{2n}{2j} \text{.}
\end{equation}
The second of the lower indices in $\alpha^{(2n)}_{2j,2(n - j)}$ is not strictly necessary, as it
is defined by the first lower and the upper indices. However, it makes for convenient notation
relating the indices of $\alpha$ with the indices of $\bm{p}$ in the outer product.

When $h \to 0$, by \cref{eq:deca:h0:Odd,eq:decb:h0:Odd} the matrix $\bm{A}$ is decreasing faster
than the right-hand side $\bm{b}$, and
therefore, the solution $\bm{w}$ must grow on the order of at least $h^{-2}$. The sum of these
increasing weights equals to one, and therefore, some of them must be negative. Thus,
no all-positive solution may exist for very small $h$. From \cref{sec:LargeH:Odd},
all weights are positive for large $h$. Therefore, for each $M$ there must be some minimal
`cutoff' value $h_\text{cut}$, such that a non-negative solution is possible only for $h \ge
h_\text{cut}$.

The weights $\bm{w}$ may grow faster than $h_*^{-2}$ if these rapidly growing terms are
orthogonal to $\bm{b}$, to make the $L^2$ finite. In \cite{Mikhin:2026:Asymptotic} we demonstrate that the
lowest power of $h_*$ is $-2M$. Then the weights are represented as:
\begin{equation}
    \label{eq:decw:h0:Odd}
        \bm{w} = \sum_{m = -M}^\infty \bm{w}_{2 m} h_*^{2 m} = \sum_{m = -M}^\infty h_*^{2 m} \sum_{j=1}^M C_{2j}^{(2m)} \tilde{\bm{p}}_{2j}  \text{.}
\end{equation}
Here we introduced complementary vectors $\tilde{\bm{p}}_{2m}$, $m = 1, ..., M$, that are orthogonal to all
$\bm{p}_{2j}$ for $j < M$ and $j \ne m$. The complementary vectors are normalised
such that $\langle \tilde{\bm{p}}_{2m}, \bm{p}_{2m}
\rangle = 1$, where angular brackets denote the inner product. The $\tilde{\bm{p}}_{2m}$ vectors are also independent, but not mutually orthogonal;
they are no longer orthogonal to $\bm{p}_{2j}$ when $j > M$.

Substituting the series \labelcref{eq:deca:h0:Odd,eq:decb:h0:Odd,eq:decw:h0:Odd} into the system
$\bm{A} \bm{w} = \bm{b}$ and equating the terms with the same powers of $h_*$, we obtain a linear
system for $C_{2j}^{(2m)}$.
Its solution is intricate because the lower-order
matrices $\bm{A}_{2n}$ are not full-rank:
the rank is defined by the order $n$ and dimensionality $M$. For
$\bm{A}_4$, all rows are multiples of the first,
and therefore, $\bm{A}_4$ has rank $1$.
Similarly, $\bm{A}_6$ has rank $2$ and $\bm{A}_8$ has rank $3$ (assuming $M \ge 3$). For arbitrary $n$, the rank of
$\bm{A}_{2n}$ is $\min{\! (n - 1, M)}$. We identify solvable subsets of
linear equations for $C_{2j}^{(2m)}$ across multiple powers of $h_*$, and obtain
\begin{align}
    \label{eq:L2h0:Odd:main}
        L^2_0 &= {\lVert F \rVert}_2 - \sum_{n=1}^M \beta_{2n} C_{2n}^{(-2n)} \text{,} \\
    \label{eq:L2h2:Odd:main}
        L^2_2 &= -\sum_{n=1}^M \beta_{2n} C_{2n}^{(-2(n - 1))} - \beta_{2(M + 1)} C_{2M}^{(-2M)} \langle \bm{p}_{2(M + 1)}, \tilde{\bm{p}}_{2M} \rangle \text{.}
\end{align}
The coefficients $C_{2n}^{(-2(n - k))}$ for $n = 1, \ldots, M$, $k = 0, 1$, are found by solving
\begin{equation}
    \label{eq:L2h0h2:Odd:C}
        \bm{G} \bm{s}_{2k} = \bm{r}_{2k} \text{,}
\end{equation}
where
\begin{equation}
    \label{eq:L2h0h2:Gsc}
    \begin{aligned}
        \bm{G} &= \Big \| \alpha^{(2(i + j))}_{2i, 2j} \Big \| \quad \text{for $i, j = 1, \ldots, M$} \text{,} \\
        \bm{s}_{2k} &= \left[ C_2^{(-2(1 - k))}, C_4^{(-2(2 - k))}, \dots, C_{2M}^{(-2(M - k))} \right] \text{,} \\
        \bm{r}_0 &= \bigl[ \beta_2, \beta_4, \dots, \beta_{2M} \bigr] \text{,} \\
        \bm{r}_2 &= \left( -C_{2M}^{(-2M)} \!
              \begin{bmatrix} \alpha^{(2(M + 2))}_{2, 2(M + 1)}, \alpha^{(2(M + 3))}_{4, 2(M + 1)}, \dots, \alpha^{(2(2M + 1))}_{2M, 2(M + 1)} \end{bmatrix} \right. + \\
              &\quad \left. \begin{bmatrix} 0, \dots, 0, \beta_{2(M + 1)} - \sum_{m=1}^M \alpha^{2(m + M + 1)}_{2(M + 1), 2m} C^{(-2m)}_{2m} \end{bmatrix} \right) \langle \bm{p}_{2(M + 1)}, \tilde{\bm{p}}_{2M} \rangle \text{.}
    \end{aligned}
\end{equation}
The coefficients $C_{2n}^{(-2n)}$ in $\bm{r}_2$ are already known from $\bm{s}_0$.
The higher-order terms $L^2_{2k}$ for $k > 1$ are given in \cite{Mikhin:2026:Asymptotic}.

The $L^2_0$ is positive, as this is the limit of $L^2$, a positive quantity, for $h \to 0$. It is
natural to expect that $L^2_2$ is negative, as having some spread of mixands must provide a better
approximation than having no spread at all. However, we currently do not have a formal proof of this
conjecture.

\subsubsection{Asymptotic solution for large \texorpdfstring{$M$}{M}}
\label{sec:LargeM:Odd}

Another asymptotic approximation can be found in the limit of $M \to \infty$. First,
we observe that\footnote{For the remainder of this derivation we use the shorter subscript form to
differentiate between the functions of different arguments.}
\begin{align}
    \notag
        \mathcal{N}_x(m h, \sigma^2) &= \frac{1}{\sqrt{2 \pi \sigma^2}}
            \exp{\! \left[ -\frac{ {(m - x / h)}^2 }{ 2 {(\sigma / h)}^2 } \right]} \\
    \label{eq:XtoM:Odd}
            &= \frac{ \mathcal{N}_m(x / h, \sigma^2 / h^2) }{h} \text{.}
\end{align}
We seek the asymptotic expression for the weights in the Gaussian form
\begin{equation}
    \label{eq:AsInfM}
    w_m \approx C \mathcal{N}_m(\gamma, \tau^2)
\end{equation}
for some $C$, $\gamma$, and $\tau$ to be determined; from the symmetry considerations, we can
set the mean $\gamma$ to zero, but keep the expression in \cref{eq:AsInfM} general for
future use. The product of the two Gaussians, now functions of $m$, is
\begin{align}
    \notag
        w_m \mathcal{N}_x(m h, \sigma^2) &\approx
            \frac{ C / h }{ \sqrt{2 \pi ( \tau^2 + \sigma^2 / h^2 ) } }
            \exp{ \left[ -\frac{ {(x / h - \gamma)}^2 }{ 2 ( \tau^2 + \sigma^2 / h^2 ) } \right] }
            \mathcal{N}_m(m_c, \sigma_c^2) \\
    \label{eq:SplitOffM:Odd}
            &= C \mathcal{N}_x \! \left( \gamma h, \sigma^2 + \tau^2 h^2 \right)
            \mathcal{N}_m(m_c, \sigma_c^2)
            \text{,}
\end{align}
where the mean and variance in $m$ are
\begin{equation}
    \label{eq:statec:Odd}
        m_c = \sigma_c^2 \left(\frac{\gamma}{\tau^2} + \frac{xh}{\sigma^2}\right)
    \quad \text{and} \quad
        \sigma_c^2 = \left(\frac{1}{\tau^2} + \frac{h^2}{\sigma^2}\right)^{-1}
    \text{.}
\end{equation}
Substituting \cref{eq:SplitOffM:Odd} into \cref{eq:Approx} and applying the Euler-Maclaurin formula
\cite[p. 16]{Abramowitz:1965:Handbook} (without the derivative terms), we replace the sum with an
integral over the same variable $m$. Being a PDF, $\mathcal{N}_m(m_c, \sigma_c^2)$ integrates to
one, and therefore, the approximation becomes
\begin{equation}
    \notag
    \widetilde{\mathcal{N}}_x(0, 1) \approx C \mathcal{N}_x \! \left( \gamma h, \sigma^2 + \tau^2 h^2 \right) \text{,}
\end{equation}
yielding
\begin{equation}
    \label{eq:AsInfM:Odd}
    C = 1 \text{,} \quad \tau = \sqrt{1 - \sigma^2} / h \text{, } \quad \gamma = 0 \text{,}
\end{equation}
and also validating our choice of the Gaussian expression for weights in \cref{eq:AsInfM}.
The expressions \labelcref{eq:statec:Odd} for the mean and variance simplify to
\begin{equation}
    \label{eq:statec:Odd:2}
        m_c = \frac{x(1 - \sigma^2)}{h}
    \quad \text{and} \quad
        \sigma_c^2 = \frac{\sigma^2 (1 - \sigma^2)}{h^2}
    \text{.}
\end{equation}
Note that $\sigma_c^2$ is a constant, while $m_c$ is a linear function of $x$.

To compute the asymptotic value of $L^2$ as $M \to \infty$, we substitute
\cref{eq:SplitOffM:Odd} into \cref{eq:L2} and take the limit:
\begin{align}
    \notag
        \lim_{M \to \infty} L^2 &= \int_{-\infty}^{+\infty}
            {\left[ \mathcal{N}_x(0, 1) -
                      C \mathcal{N}_x \! \left(\gamma h, \sigma^2 + \tau^2 h^2 \right) \lim_{M \to \infty} \sum_{m = -M}^{M}
            \mathcal{N}_m(m_c, \sigma_c^2) \right]}^2 \! dx \\
    \label{eq:AsympL2:Odd}
                                &= \int_{-\infty}^{+\infty} \mathcal{N}_x(0, 1)^2
            {\left[ 1 - \sum_{m = -\infty}^{+\infty} \mathcal{N}_m(m_c, \sigma_c^2) \right]}^2 \! dx
            \text{.}
\end{align}
To find the infinite sum,
consider the following result from \cite[p.~481]{Graham:2017:Concrete}:
\begin{equation}
    \label{eq:AsympExp}
    \sum_{k=-\infty}^{+\infty} e^{-\frac{(k + t)^2}{n}} = \sqrt{\pi n}\left(1 + 2 \sum_{m=1}^{+\infty} e^{-\pi^2 n m^2} \cos{(2 \pi m t)}\right)\text{.}
\end{equation}
The left-hand side of \cref{eq:AsympExp} matches $\mathcal{N}_m(m_c, \sigma_c^2)$ for $k
= m$, $t = -m_c$ and $n = 2 \sigma_c^2$.
The exponential term in each element decays rapidly as $m \to \infty$, while the cosine term is
bounded on $[-1, 1]$. Thus, we approximate the series by its first element as:
 \begin{equation}
    \notag
    \sum_{m = -\infty}^{+\infty}
            \mathcal{N}_m(m_c, \sigma_c^2) \approx 1 + 2 e^{-2 \pi^2 \sigma_c^2} \cos{(2 \pi m_c)} + \mathcal{O}(e^{-8 \pi^2 \sigma_c^2})\text{.}
\end{equation}
Substituting this expression into \cref{eq:AsympL2:Odd} yields, omitting the Big-$\mathcal{O}$ terms:
\begin{equation}
    \label{eq:L2:Odd:AsympM:1}
        \lim_{M \to \infty} L^2 = 4 e^{-4 \pi^2 \sigma_c^2}
            \int_{-\infty}^{+\infty} \mathcal{N}_x(0, 1)^2 \cos^2{\! (2 \pi m_c)} \, dx \text{.}
\end{equation}
Using the cosine double angle identity and applying
\cite[Eq.~7.4.6]{Abramowitz:1965:Handbook}, we take this integral to find
the
asymptotic value of $L^2$ as:
\begin{equation}
    \label{eq:L2:Odd:AsympM}
        \lim_{M \to \infty} L^2 \approx
            \frac{1}{\sqrt{\pi}} e^{-4 \pi^2 \frac{\sigma^2 (1 - \sigma^2)}{h^2}} \left( 1 + e^{-4 \pi^2 \frac{(1 - \sigma^2)^2}{h^2}} \right) \text{.}
\end{equation}
More details on this derivation are given in \cite{Mikhin:2026:Asymptotic,Xiourouppa:2027:PhD}.

\cref{eq:L2:Odd:AsympM} suggests that for $h \to 0$ the mismatch $L^2$ tends to zero. However, from
\cref{eq:L2h0:Odd:main} we know that the actual limit is finite. Therefore, \cref{eq:L2:Odd:AsympM}
is not applicable for very small $h$.
Still, if we assume that the areas of validity for the small-$h$ and large-$M$ approximations
overlap, we can combine these results to find an approximation for $h_\text{opt}$. Following
\cref{sec:SmallH:Odd}, at small $h$ we represent $L^2$ as the second-order power series $L^2_0 +
L^2_2 h^2$, which is expected to decrease with $h$. On the other hand, the $L^2$ approximation by
\cref{eq:L2:Odd:AsympM} tends to zero when $h \to 0$, and then grows with $h$.
Therefore, as $h$ increases the two functions must intersect, and their intersection would yield an
approximation to $h_\text{opt}$.

The solution is not possible in closed form, but can be approximated. Given that the large-$M$
prediction of $L^2$ is exponentially small, we can replace it with zero. Then the hypothetical point
of optimum becomes
\begin{equation}
    \label{eq:hdagger:Odd}
    h^\dagger_\text{opt} = - {( L^2_0 / L^2_2 )}^{\sfrac{1}{2}} \text{.}
\end{equation}
Having $h^\dagger_\text{opt}$, we then obtain the corresponding weights and $L^2$ from the large-$M$
asymptotic expressions.

\subsubsection{Heuristic solution for large \texorpdfstring{$M$}{M}, small \texorpdfstring{$h$}{h}}
\label{sec:Heur:Odd}

The large-$M$ approximation obtained in \cref{sec:LargeM:Odd} is not applicable for
small $h$.
Looking at the derivation details, the solution
requires not only $M$, but the product $M h$ to be large, such that we can
approximate the infinite integral with a sum. This more restrictive assumption fails if $h \to 0$.

To specifically address the case of intermediate $M h$, consider the following heuristic approach.
The highest-order mixand has the mean $M h$. Having many mixands from $-M h$ to $+M h$, we can
expect good approximation in that interval, so assume that its contribution to $L^2$ is zero. On the
other hand, for $|x| > M h$ we only have the tails of the mixands that, for $\sigma < 1$, decay
faster than the original Gaussian. One option is to assume that these tails provide no approximation at all, and
the value of the Gaussian sum is zero. Under such assumption, we can express the overall
$L^2$ mismatch through $\erf$, the Gauss error function \cite[Sect. 7.1]{Abramowitz:1965:Handbook}.
Potentially more accurately, we can assume that the
Gaussian sum beyond $M h$ is represented by the single highest-order component with the weight $w_M$.
This way, we obtain
\begin{align}
    \notag
        L^2_\text{heur}
        &= 2 \int_{M h}^{+\infty} \left[
            \widetilde{\mathcal{N}} - \frac{w_M}{\sqrt{2 \pi \sigma^2}} \exp{\left( -\frac{ (x - M h)^2 }{ 2 \sigma^2 } \right)}
        \right]^2 dx \\
    \label{eq:Heur:Odd}
        &= \frac{1}{2 \sqrt{\pi}} \left[ 1 - \erf{(M h)} \right]
           + \frac{w_M^2}{2 \sigma \sqrt{\pi}}
           - \frac{w_M \sqrt{2}}{\sqrt{\pi (1 + \sigma^2)}} e^{-\eta^2}
                     \left[ 1 - \erf{\left( \sigma \eta \right)} \right]
        \text{,}
\end{align}
where $\eta = M h / \sqrt{2 (1 + \sigma^2)}$. Then we select $w_M$ to minimise the error
of this approximation, yielding:
\begin{equation}
    \notag
    w_M = \frac{\sigma \sqrt{2}}{\sqrt{1 + \sigma^2}} e^{-\eta^2} \left[ 1 - \erf{\left( \sigma \eta \right)} \right]
    \text{.}
\end{equation}
With this method, we can find an approximation $h^*_\text{opt}$ to the exact optimal solution
$h_\text{opt}$ as the intersection of \cref{eq:Heur:Odd} with
the large-$M$ expression for $L^2$ from \cref{eq:L2:Odd:AsympM}. However, there is no closed form
solution.

We may further simplify the error function via the approximation
\begin{equation}
    \label{eq:erfc}
    1 - \erf{(x)} = \erfc{(x)} \approx \frac{e^{-x^2}}{x \sqrt{\pi}} \text{,}
\end{equation}
and therefore,
\begin{equation}
    \label{eq:Heur:wMAlt:Odd}
    w_M \approx \frac{\sqrt{2} e^{-\eta^2 (1 + \sigma^2)}}{\eta \sqrt{\pi (1 + \sigma^2)}}
        = \frac{2}{M h \sqrt{\pi}} e^{ -M^2 h^2 / 2 }
    \text{,}
\end{equation}
and
\begin{equation}
    \label{eq:Heur:Alt:Odd}
    \begin{aligned}
        L^2_\text{heur}
        &\approx \frac{e^{-M^2 h^2}}{2 M h \pi}
            \left( 1 - \frac{4}{\sigma M h \sqrt{\pi}} \right) \text{.}
    \end{aligned}
\end{equation}
Regrettably, this still does not allow a closed-form solution for $h^*_\text{opt}$.

\subsection{Even case}
\label{sec:1D:Even}

Now consider the complementary case where the number of terms in the mixture is even.
Again, it is convenient to change the notation as $M \Rightarrow 2 M$. Compared to the odd case, it is still more
challenging to find the indexing scheme that exposes the symmetry of components with corresponding
positive and negative mean values. Writing the approximation as
\begin{equation}
    \label{eq:Approx:Even}
    \widetilde{\mathcal{N}}(x; 0, 1) \approx \sum_{m = -M}^{M - 1} w_{m} \mathcal{N}(x; (m + \sfrac{1}{2}) h, \sigma^2) \text{,}
\end{equation}
the fully-symmetric form would need to use non-integer indices $m + \sfrac{1}{2}$ in place of $m$.
The terms $m - 1$ and $-m$ in \cref{eq:Approx:Even} have symmetric offsets from zero equal to $(m -
\sfrac{1}{2}) h$ and $(-m + \sfrac{1}{2}) h$, and therefore, must have the same weights. The
weight normalisation condition \labelcref{eq:WNorm} yields
\begin{equation}
    \notag
    \sum_{m = -M}^{M - 1} w_{m} = 2 \sum_{m = 0}^{M - 1} w_{m} = 2 w_0 + 2 \sum_{m = 1}^{M - 1} w_{m} = 1 \text{.}
\end{equation}
Using these relationships, we eliminate all non-positive indices from \cref{eq:Approx:Even} to
obtain
\begin{equation}
    \label{eq:L2:Even}
    L^2 = \int_{-\infty}^{+\infty}
        {\left[ \widetilde{\mathcal{N}} - \frac{1}{2} \mathcal{N}_{\pm \sfrac{1}{2}} -
                \sum_{m = 1}^{M - 1} w_m \left( \mathcal{N}_{-m - \sfrac{1}{2}} + \mathcal{N}_{m + \sfrac{1}{2}}
                                                -\mathcal{N}_{\pm \sfrac{1}{2}} \right) \right]}^2 \! dx \text{,}
\end{equation}
where we introduced $\mathcal{N}_{\pm \sfrac{1}{2}} \defeq \mathcal{N}_{-\sfrac{1}{2}} +
\mathcal{N}_{\sfrac{1}{2}}$. Using a new set of basis
functions
\begin{equation}
    \label{eq:Base:Even}
    f_m \defeq \mathcal{N}_{-m - \sfrac{1}{2}} + \mathcal{N}_{m + \sfrac{1}{2}} - \mathcal{N}_{\pm \sfrac{1}{2}} \text{,}
\end{equation}
and defining
\begin{equation}
    \label{eq:RHS:Even}
    F \defeq \widetilde{\mathcal{N}} - \frac{1}{2} \mathcal{N}_{\pm \sfrac{1}{2}} \text{,}
\end{equation}
we again express the $L^2$ mismatch in the form of \cref{eq:L2:Quad}, with the coefficients now
given by
\begin{equation}
    \label{eq:cd:Even}
    \begin{aligned}
        {\lVert F \rVert}_2 &= \sfrac{1}{2 \sqrt{\pi}} - 2 d_{\sfrac{1}{2}}
            + \sfrac{1}{2} \bigl( c_{\sfrac{1}{2},-\sfrac{1}{2}} + c_{\sfrac{1}{2},\sfrac{1}{2}} \bigr) \text{,}
            \\
        b_m &= 2 \bigl( d_{m + \sfrac{1}{2}} - d_{\sfrac{1}{2}} \bigr)
            - c_{-\sfrac{1}{2},m + \sfrac{1}{2}} - c_{\sfrac{1}{2},m + \sfrac{1}{2}}
            + c_{-\sfrac{1}{2},\sfrac{1}{2}} + c_{\sfrac{1}{2},\sfrac{1}{2}}
            \text{,}
            \\
        a_{m,k} = & +2 \bigl( c_{m + \sfrac{1}{2},k + \sfrac{1}{2}} + c_{m + \sfrac{1}{2},-k - \sfrac{1}{2}} \bigr)
        \\
                 & -2 \bigl( c_{-\sfrac{1}{2},k + \sfrac{1}{2}} + c_{-\sfrac{1}{2},m + \sfrac{1}{2}} \bigr)
        \\
                 & -2 \bigl( c_{\sfrac{1}{2},k + \sfrac{1}{2}} + c_{\sfrac{1}{2},m + \sfrac{1}{2}} \bigr)
        \\
                 & +2 \bigl( c_{-\sfrac{1}{2},\sfrac{1}{2}} + c_{\sfrac{1}{2},\sfrac{1}{2}} \bigr)
            \text{.}
            \\
    \end{aligned}
\end{equation}
We heavily used the symmetry rules for $d_\alpha$ and $c_{\alpha,\beta}$, but the resulting
expressions are still notably more onerous than for the odd case in \cref{eq:cd:Odd}. Solving the obtained system of
linear equations, we find $\bm{w}(h) = \mathbf{A}^{-1} \bm{b}$, and obtain $L^2$ as a
univariate function of $h$. This solution does not cover the special case of $M =
1$, \textit{i.e.}, a two-term approximation. However, the weights for this corner case are
immediately available from the symmetry of approximation: $w_{-1} = w_0 = \sfrac{1}{2}$. Then the
mismatch norm equals ${\lVert F \rVert}_2$ from \cref{eq:cd:Even}.

\subsubsection{Asymptotic solution for large \texorpdfstring{$h$}{h}}
\label{sec:LargeH:Even}

To analyse the asymptotic behaviour of the obtained solution, recall that the lower subscripts
$\alpha$, $\beta$ in $d_\alpha$, $c_{\alpha,\beta}$ are not just indices but
\textit{multipliers} of $h^2$ per \cref{eq:Gprods}. For the odd case, the multipliers were integers,
but not here. Thus, in the limit of $h \to \infty$, we have $c_{\alpha,\beta} \to 0$ unless $\alpha
= \beta$, $d_\alpha \to 0$ unless $\alpha = 0$, and the terms in \cref{eq:cd:Even} simplify to
\begin{equation}
    \notag
    \begin{aligned}
        b_m &\to c_{\sfrac{1}{2},\sfrac{1}{2}}
            &&= \frac{1}{2 \sigma \sqrt{\pi}} \text{,} \\
        a_{m,k} &\to 2 c_{m + \sfrac{1}{2},k + \sfrac{1}{2}} \delta_{m,k} + 2 c_{\sfrac{1}{2},\sfrac{1}{2}}
            &&= \frac{1 + \delta_{m,k}}{\sigma \sqrt{\pi}} \text{,}
    \end{aligned}
\end{equation}
or
\begin{equation}
    \notag
    \mathbf{A} = \frac{1}{\sigma \sqrt{\pi}} (\mathbf{J} + \mathbf{I}) \quad \text{and} \quad
    \bm{b} = \frac{1}{2 \sigma \sqrt{\pi}} \bm{j} \text{.}
\end{equation}
From the normalisation condition, we have $\mathbf{J} \bm{w} = (1/2 - w_0) \bm{j}$, and therefore,
the equation for $\bm{w}$ reduces to $\mathbf{I} \bm{w} = w_0 \bm{j}$. Thus, all weights are the
same, $w_m = \sfrac{1}{2 M}$, $\forall m$. We also obtain the $L^2$ norm as
\begin{equation}
    \label{eq:AsInf:Even}
    L^2_\infty = \frac{1}{2 \sqrt{\pi}} + \frac{1}{4 M \sigma \sqrt{\pi}} \text{.}
\end{equation}
An intuitive explanation for this expression is that in the limit of $h \to
\infty$, all the even-case mixands ``move'' to infinity and do not overlap with the original
distribution. The squared $L^2$ norm of the mismatch is then the squared norm of the original
distribution itself plus the squared norm of a mixand multiplied by its weight (squared) and taken
$2 M$ times, yielding \cref{eq:AsInf:Even}.

\subsubsection{Asymptotic solution for small \texorpdfstring{$h$}{h}}
\label{sec:SmallH:Even}

Due to the proliferation of terms in \cref{eq:cd:Even} in comparison to \cref{eq:cd:Odd}, the asymptotic analysis of the even-length
case for $h \to 0$ is much more cumbersome. The details are given in
\cite{Mikhin:2026:Asymptotic}. We again search for the solution in the form of power
series over $h_*$, as in \cref{sec:SmallH:Odd}. For the right-hand side
vector $\bm{b}$ we now obtain (\textit{cf.} \cref{eq:decb:h0:Odd,eq:beta:h0:Odd,eq:p2k:Odd})
\begin{equation}
    \label{eq:decb:h0:Even}
    \bm{b} = \sum_{n=1}^{+\infty} h_*^{2n} \sum_{j = 1}^{n} \beta^{(2j)}_{2n} \bm{p}_{2j} \text{,}
\end{equation}
where
\begin{equation}
    \label{eq:beta:h0:Even:2n}
         \beta^{(2n)}_{2n} = \frac{(-1)^n}{n! \, 2^{2n} \sigma \sqrt{\pi}} \left( \frac{2^{n + 1/2} \sigma^{2n + 1}}{(1 + \sigma^2)^{n + 1/2}} - 1 \right)
         \text{,}
\end{equation}
which is the same as the odd-case value from \cref{eq:beta:h0:Odd},
\begin{equation}
    \label{eq:beta:h0:Even:2j}
         \beta^{(2j)}_{2n} = \frac{(-1)^{n + 1}}{n! \, 2^{4n - 2j} \sigma \sqrt{\pi}} \binom{2n}{2j} \text{,}
         \quad \forall j < n
         \text{,}
\end{equation}
and the vectors $\bm{p}_{2k}$ are now defined as (\textit{cf.} \cref{eq:p2k:Odd})
\begin{equation}
    \label{eq:p2k:Even}
    \bm{p}_{2k} = \bigl[
        (1 + \sfrac{1}{2})^{2k} - (\sfrac{1}{2})^{2k},
        \dots,
        (M + \sfrac{1}{2})^{2k} - (\sfrac{1}{2})^{2k}
    \bigr] \text{, } \forall k \ge 1 \text{.}
\end{equation}
For the left-hand matrix $\bm{A}$ we again obtain the power-series decomposition
by \crefrange{eq:deca:h0:Odd}{eq:alpha:h0:Odd}, just with the new definition of the $\bm{p}_{2k}$
vectors from \cref{eq:p2k:Even} instead of \cref{eq:p2k:Odd}.

As in the odd case, the power series of $\bm{b}$ starts at $h_*^2$, while the series of $\bm{A}$
only at $h_*^4$. Therefore, the mixand weights grow infinitely by modulo as $h \to 0$, some of
them become negative to comply with the normalisation condition, and for each $M$ there is
some minimal `cutoff' value $h_\text{cut}$, such that a non-negative solution is possible only for
$h \ge h_\text{cut}$.

We search for the weights in the power series form per \cref{eq:decw:h0:Odd}, with the
complementary vectors $\tilde{\bm{p}}_{2m}$, $m = 1, ..., M$, defined in the same way as before,
just using the new vectors $\bm{p}_{2k}$ from \cref{eq:p2k:Even}. The solution is then fundamentally
equivalent to the odd case: equating the terms with the same powers of $h_*$ in $\bm{A} \bm{w} =
\bm{b}$, we obtain a linear system of equations for the coefficients $C_{2j}^{(2m)}$, and
incrementally solve it. The explicit form of the solution is more unwieldy than in
\cref{eq:L2h0h2:Gsc}. Due to the presence of multiple vectors $\bm{p}_{2j}$, $j \ne n$, in
\cref{eq:decb:h0:Even}, the obtained equations have more non-zero terms, resulting in
additional nested sums in $\bm{r}_0$, $\bm{r}_2$ of the even-case equivalent of
\cref{eq:L2h0h2:Gsc}. Thus, we omit these equations here and refer to
\cite{Mikhin:2026:Asymptotic} for the details.

\subsubsection{Asymptotic solution for large \texorpdfstring{$M$}{M}}
\label{sec:LargeM:Even}

We again look for the asymptotic approximation for the weights in the limit of $M \to \infty$. In
place of \cref{eq:XtoM:Odd}, we find
\begin{equation}
    \label{eq:XtoM:Even}
        \mathcal{N}_x((m + \sfrac{1}{2}) h, \sigma^2)
            = \frac{ \mathcal{N}_m(x / h - \sfrac{1}{2}, \sigma^2 / h^2) }{h} \text{.}
\end{equation}
Using the asymptotic form for the weights from \cref{eq:AsInfM} and transforming the product
of two Gaussians as in \cref{eq:SplitOffM:Odd}, we get
\begin{equation}
    \label{eq:SplitOffM:Even}
    w_m \mathcal{N}_x(m h, \sigma^2) \approx
        C \mathcal{N}_x \! \left( (\gamma + \sfrac{1}{2}) h, \sigma^2 + \tau^2 h^2 \right)
        \mathcal{N}_m(m_c, \sigma_c^2)
        \text{.}
\end{equation}
By substituting this expression into \cref{eq:Approx} and integrating over $m$ using the
Euler-Maclaurin formula \cite[p.~16]{Abramowitz:1965:Handbook}, we obtain the approximation
parameters
\begin{equation}
    \label{eq:AsInfM:Even}
    C = 1 \text{,} \quad \tau = \sqrt{1 - \sigma^2} / h \text{, } \quad \gamma = -1 / 2 \text{.}
\end{equation}
The values of $C$ and $\tau$ are the same as in the odd case (\cref{eq:AsInfM:Odd}); only $\gamma$, the
mean, differs, which is again a consequence of the asymmetry in the integer-valued indices: for $m =
-M, \ldots, M - 1$, we have $m - \gamma = m + \sfrac{1}{2}$ iterating over the symmetric range of
$-M + \sfrac{1}{2}, \ldots, M - \sfrac{1}{2}$.

The derivation of $L^2$ as $M \to \infty$ is similar to the odd case, but there are two minor
differences. First, the upper limit of the sum in \cref{eq:Approx} is $M-1$ instead of $M$, which
has no impact since $M-1 \to \infty$ as $M \to \infty$. Second, the value of $m_c$ becomes
\begin{equation}
    \label{eq:mc:Even}
        m_c = \frac{x(1 - \sigma^2)}{h} - \frac{1}{2} \text{.}
\end{equation}
As earlier, it is a linear function of $x$; $\sigma_c$ does not change.
Substituting \cref{eq:SplitOffM:Even} into \cref{eq:Approx} and applying the Euler-Maclaurin
formula,
we find the $L^2$ error as
\begin{align}
    \notag
        \lim_{M \to \infty} \! L^2 &= \int_{-\infty}^{+\infty}
            \! {\biggl[ \mathcal{N}_x(0, 1) -
                \! \sum_{m = -\infty}^{+\infty} \! C \mathcal{N}_x \! \left( (\gamma + \sfrac{1}{2}) h, \sigma^2 + \tau^2 h^2 \right)
            \mathcal{N}_m(m_c, \sigma_c^2) \biggr]}^2 \! dx \\
    \notag
            &= \int_{-\infty}^{+\infty} \mathcal{N}_x \! \left(0, 1\right)^2
            {\left[ 1 - \sum_{m = -\infty}^{\infty} \! \mathcal{N}_m(m_c, \sigma_c^2) \right]}^2
            \! dx \\
    \label{eq:L2:Even:AsympM:1}
            &\approx 4 e^{-4 \pi^2 \sigma_c^2}
            \int_{-\infty}^{+\infty} \mathcal{N}_x(0, 1)^2 \cos^2{\! (2 \pi m_c)} \, dx \text{.}
\end{align}
This expression is identical to \cref{eq:L2:Odd:AsympM:1}, only the $m_c$ is different. We
reuse most of the derivation presented for the odd case; details are given in
\cite{Mikhin:2026:Asymptotic,Xiourouppa:2027:PhD}. The resulting asymptotic value of $L^2$ is
the same as \cref{eq:L2:Odd:AsympM}. We conclude that up to the next order terms, the
asymptotic $L^2$ errors for $M \to \infty$ are equal for odd and even $M$, which is intuitively
expected.

\subsubsection{Heuristic solution for moderate \texorpdfstring{$M h$}{Mh}}
\label{sec:Heur:Even}

The heuristic solution presented in \cref{sec:Heur:Odd} easily extends to the even case.
Now the highest-order mixand has the mean $(M - \sfrac{1}{2}) h$. For large $M h$, the
factor of $\sfrac{1}{2}$ is negligible.
There are many mixands from zero to $(M - \sfrac{1}{2}) h$, so assume that the contribution of this
interval to $L^2$
is zero, while for $|x| > (M - \sfrac{1}{2}) h$
assume that the Gaussian sum is represented by the
single highest-order component with the weight $w_{M - \sfrac{1}{2}}$. This yields
\begin{equation}
    \label{eq:Heur:Even}
        L^2_\text{heur}
        = 2 \int_{(M - \sfrac{1}{2}) h}^{+\infty} \left[
            \widetilde{\mathcal{N}} - \frac{w_{M - \sfrac{1}{2}}}{\sqrt{2 \pi \sigma^2}} \exp{\left( -\frac{ (x - (M - \sfrac{1}{2}) h)^2 }{ 2 \sigma^2 } \right)}
        \right]^2 dx
        \text{.}
\end{equation}
Comparing \cref{eq:Heur:Even} with \cref{eq:Heur:Odd}, we
observe that the only difference is the substitution of $M$ with $M - \sfrac{1}{2}$. Therefore, we
can reuse the results of \cref{sec:Heur:Odd} by applying this replacement in the final
\cref{eq:Heur:Odd,eq:Heur:Alt:Odd}. As in the odd case, we
can find an intersection of this heuristic solution with the large-$M$ asymptotic from
\cref{eq:L2:Odd:AsympM}, and use the obtained value $h^*_\text{opt}$ as an approximation of the
exact optimum $h_\text{opt}$.

\section{Numerical results}
\label{sec:Numerics:1D}

\subsection{Computational aspects}
\label{sec:Computing:1D}

The implementation of the proposed splitting algorithm consists of two distinct stages. The first
part solves the linear system for the weights $\bm{w}$.
This can be done for arbitrary step $h$ between mixand means. Then, in the second stage, we find the
minimum of the $L^2$ mismatch with respect to $h$ using the weights produced by the first part of the algorithm.
Minimisation was implemented using routines from the {\tt scipy} package \cite{Virtanen:2020:SciPy}
with no issues, whereas calculations at the first stage
exhibited numerical instability for small $h$, and insufficient precision for large $M$ and the resulting small
$L^2$.

Theoretical analysis of $L^2$ for $h \to 0$ in \cref{sec:SmallH:Odd,sec:SmallH:Even} revealed the extreme
degeneracy of the matrix $\bm{A}$ as the root cause of the first problem. The issue was
bypassed by computing and solving the linear system using high-precision floating-point arithmetic
provided by the {\tt mpmath} package \cite{mpmath}. Precision to 256 decimal places was sufficient
for all examples presented below. Application of {\tt mpmath} also resolved the
precision issues for large $M$.

The use of high-precision arithmetic was limited to the first stage only. Once the weights and the
mismatch were obtained in high-precision, the
subsequent minimisation of $L^2$ and other calculations below were performed in standard
64-bit floating point arithmetic.

\subsection{Examining the solution}
\label{sec:Exam:1D}

\begin{figure}[htb!]
     \centering
        \includegraphics[width=\textwidth]{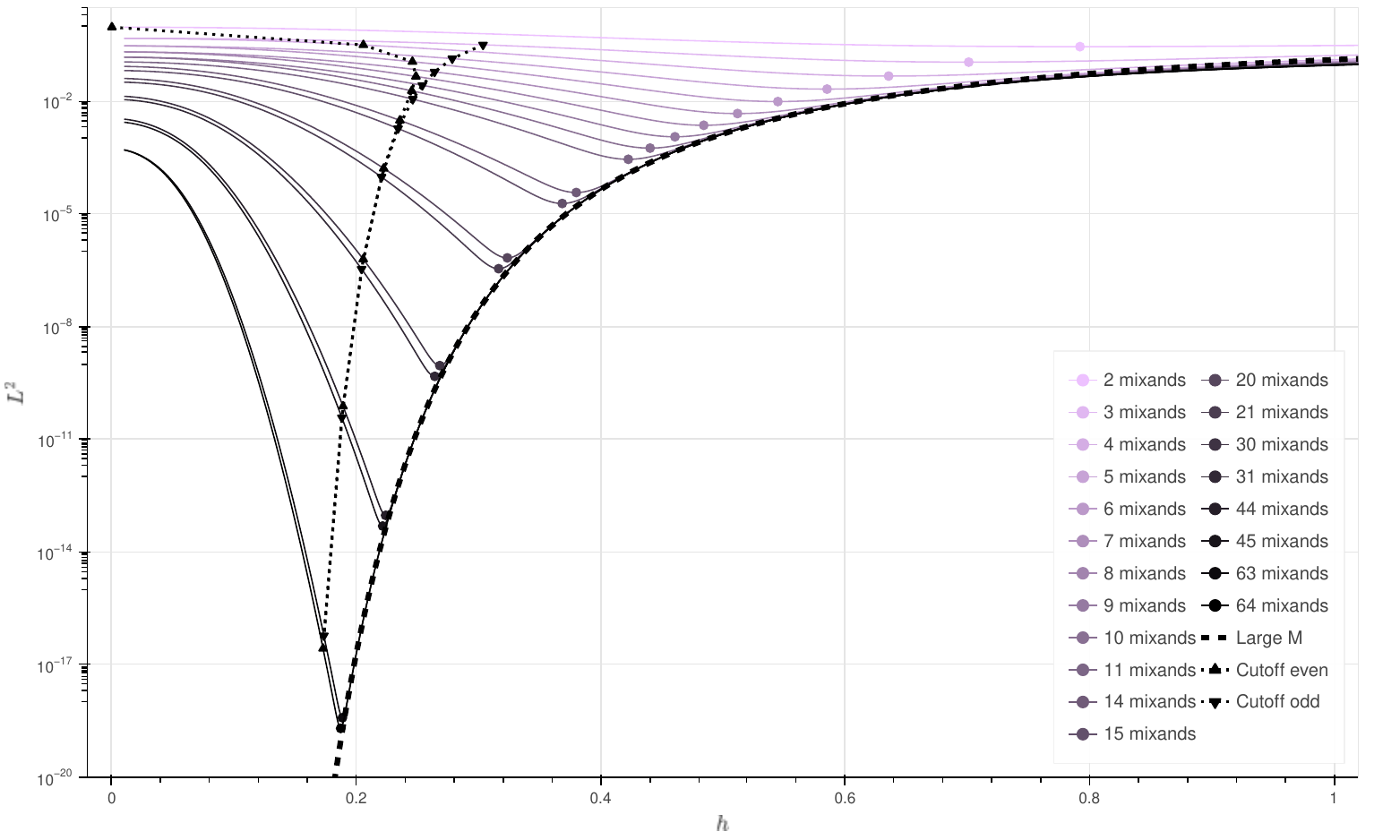}
        \caption{The approximation mismatch $L^2$ as a function of the mean step $h$, part 1. Solid
        lines: exact solutions for various $M$. Dots: points of optima. Thick dashed line:
        large-$M$ asymptotic solution. Dotted lines: cut-off boundaries for odd and even cases;
        there are no all-positive $\bm{w}$ solutions for smaller values of $h$.}
        \label{fig:L2h_M_c}
    \vspace{-3mm}
\end{figure}

\begin{figure}[htb!]
    \centering
    \includegraphics[width=\textwidth]{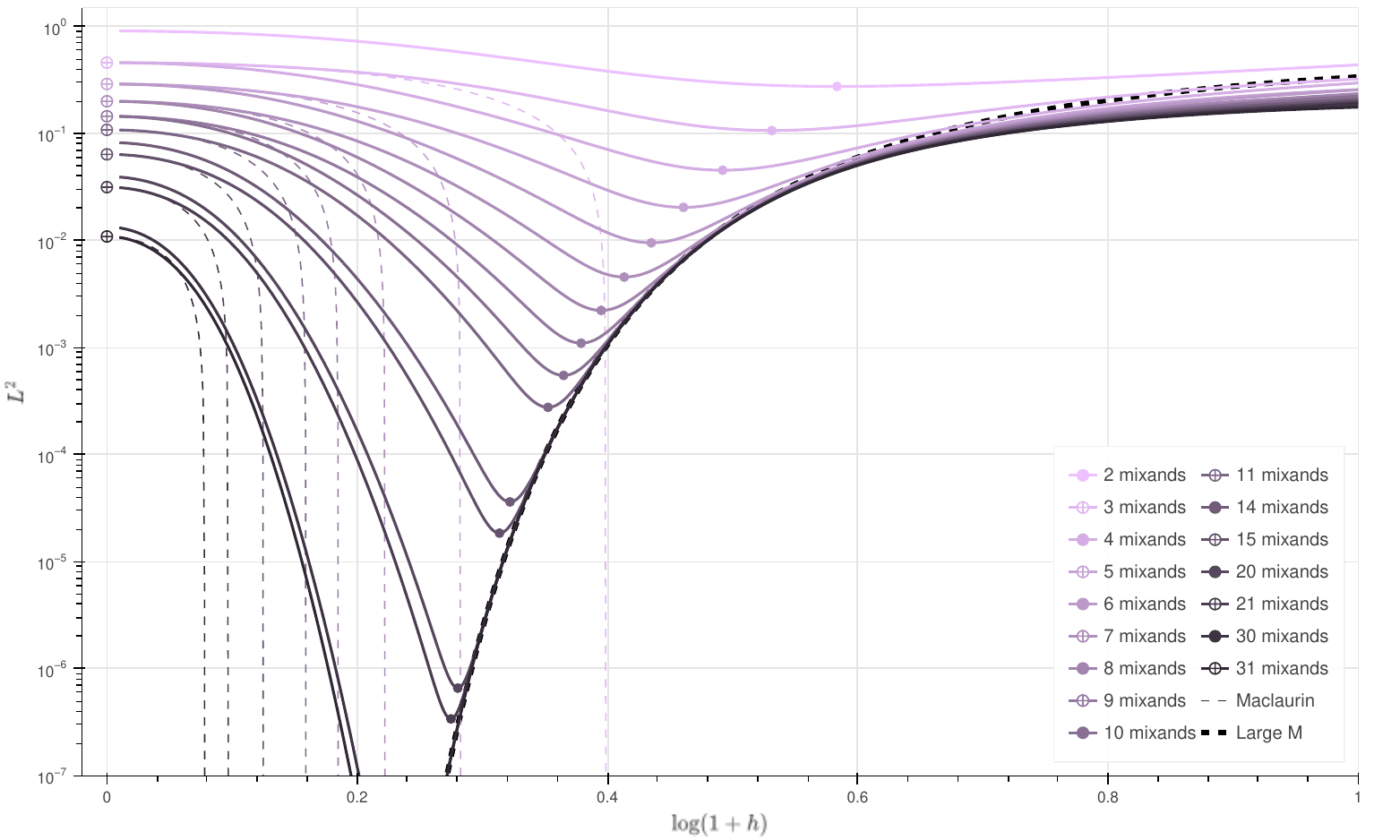}
    \caption{The approximation mismatch $L^2$ as a function of the mean step $h$, part 2. Solid
    lines, dots, and thick dashed line: as in \cref{fig:L2h_M_c}. Thin dashed lines: small-$h$
    solutions from \cref{sec:SmallH:Odd}, odd case only.}
    \label{fig:L2h_M_h0}
    \vspace{-3mm}
\end{figure}

\begin{figure}[htb!]
    \centering
    \includegraphics[width=\textwidth]{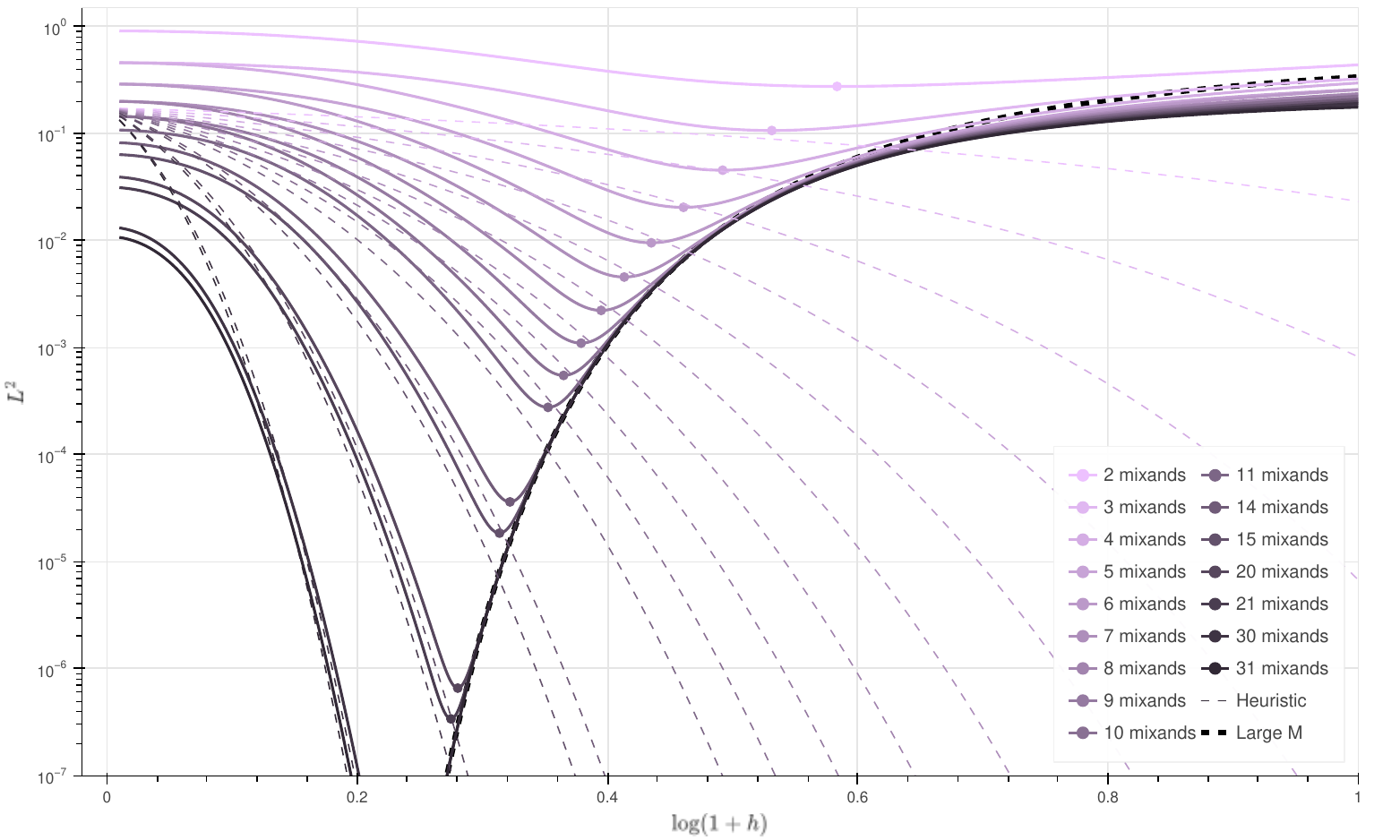}
    \caption{The approximation mismatch $L^2$ as a function of the mean step $h$, part 3. Solid
    lines, dots, and thick dashed line: as in \cref{fig:L2h_M_c}.
    Thin dashed lines: heuristic solutions from
    \cref{sec:Heur:Odd,sec:Heur:Even}.}
    \label{fig:L2h_M_hr}
    \vspace{-3mm}
\end{figure}

The core of our approach is minimising the $L^2$ mismatch between the original Gaussian and its
approximation. After solving for the weights, the mismatch is described by the function
$\mathcal{M}(h; \sigma, M)$ from \cref{eq:L2:special}. We start by examining the behaviour of this
function on $h$ and $M$
illustrated in \cref{fig:L2h_M_c,fig:L2h_M_h0,fig:L2h_M_hr}. Throughout this section, the
mixand standard deviation is $\sigma = 0.2$ unless specified otherwise.

\cref{fig:L2h_M_c} demonstrates that for all considered values of $M$, the functions $L^2(h)$ are
unimodal with unique minima $h_\text{opt}$ (dots). The mismatch values at the
optimal steps are small and steadily decrease with $M$. As discussed in \cref{sec:SmallH:Odd,sec:SmallH:Even}, the
approximations are not acceptable for $h$ values below the cutoffs $h_\text{cut}$ (dotted lines with
triangular markers), as at least some of the weights are negative; cutoff computation is discussed
below. The cutoff lines differ for odd and even $M$, but converge to the same
curve as $M$ grows. All optimal solutions are to the right of the cutoff lines, and therefore,
yield valid splittings. For large $M$, the asymptotic solution derived in
\cref{sec:LargeM:Odd,sec:LargeM:Even} (thick dashed line) provides an accurate approximation
of the exact solutions for $h \gtrsim h_\text{opt}$. Currently we do not have
corresponding asymptotic solutions for the cutoff lines. At the same time, the large-$M$ approximation
is ``one-sided'': it predicts ever decreasing values of $L^2$ as $h$ gets smaller, while the exact
solution has a minimum at the optimal $h_\text{opt}$ and then increases again.

In \cref{sec:LargeM:Odd} we discussed using
$h^\dagger_\text{opt}$ as an approximation of the true $h_\text{opt}$. The analysis assumed that
the areas of validity for the small-$h$ and large-$M$ approximations overlap. This assumption
apparently fails.
Dashed lines in \cref{fig:L2h_M_h0} show the $L^2$ approximation by two-term Maclaurin series for small $h$;
they lose accuracy and drop to zero well to the left from
$h_\text{opt}$, invalidating the idea that the small-$h$ and large-$M$ approximations can be combined to
estimate $h_\text{opt}$.
Instead, such an estimate is apparently achieved using the heuristic solutions $L^2_\text{heur}(h)$
(dashed lines in \cref{fig:L2h_M_hr})
that intersect the large-$M$ curve very close to the locations of the true minima. Despite
the empirical nature of this approximation, it is accurate for moderate $M
h$. The $L^2$ values at the intersection points are somewhat optimistic, but are usually
within one decimal order from the exact solutions.

For small $M$, we can compare our approximations to those from \cite{Faubel:2010:Further},
which are simpler, but only modestly reduce the mixand variance, which in part motivated our work.
\cref{tab:faubel_comp} shows the $L^2$ mismatches using the $\nu$ ($h$ in our notation) and $\sigma$
values recommended by \cite{Faubel:2010:Further}, as well as our optimal results for $M = 2, 3$ and
the same mixand variances.
The two methods attain equal orders of magnitude of $L^2$ for the presented parameters,
although our mismatch is, by design, smaller. Overall, the approximation by
\cite{Faubel:2010:Further} is close to optimal, and its modest reduction in the mixand variance is
just an unavoidable consequence of using so few mixands. Achieving both good approximation and
substantial reduction in variance requires more components.

\begin{table}[htb!]
    \caption{$L^2$ mismatch using our approximation and the method of \cite{Faubel:2010:Further} for the same
    mixand variances and the number of components.}
    \label{tab:faubel_comp}
    \centering
\begin{tabular}{|c|c|c|c|c|c|}
\hline
\textbf{M} & \textbf{Faubel $\nu$} & \textbf{Faubel $\sigma$} & \textbf{$L^2$ for $\nu$} & \textbf{Best $h$} & \textbf{Minimal $L^2$}  \\ \hline
\textbf{2}               & $0.5$                 & $0.8660$                 & $9.700 \times 10^{-5}$   & $0.4698$          & $2.546 \times 10^{-5}$  \\ \hline
\textbf{3}               & $0.5$                 & $0.9574$                 & $9.322 \times 10^{-10}$  & $0.4899$          & $1.118 \times 10^{-10}$ \\ \hline
\end{tabular}
\end{table}

\begin{figure}[htb!]
    \centering
    \includegraphics[width=\textwidth]{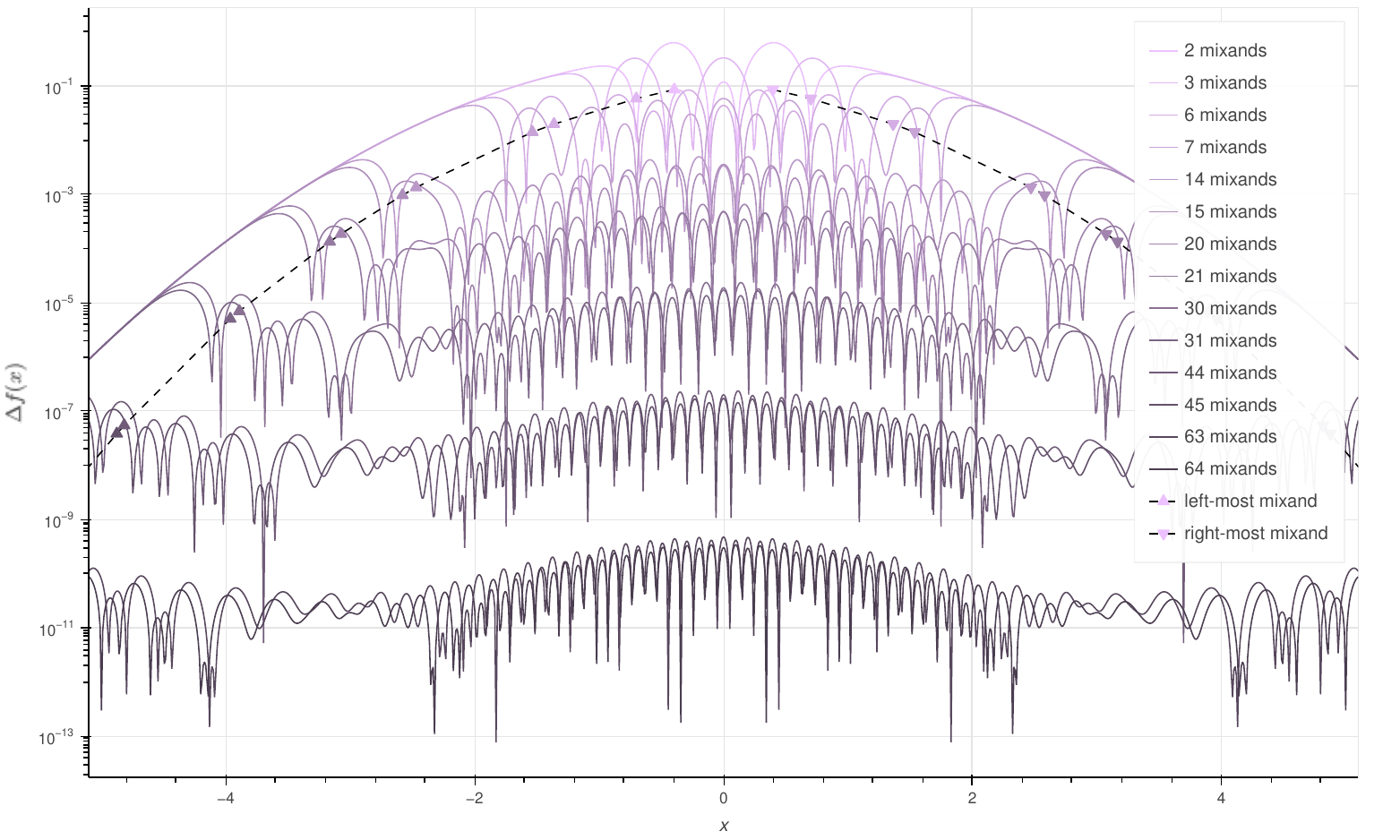}
    \caption{Solid lines: errors in approximating the probability density function of the original
    Gaussian distribution. Dashed lines with markers: positions of the outer-most mixand means
    ($x$ values) vs. the mean approximation errors ($y$ values).}
    \label{fig:GS_mism}
    \vspace{-3mm}
\end{figure}

\cref{fig:GS_mism} displays the difference between our GS approximations
$\widetilde{\mathcal{Q}}(x)$ and the original distribution $\widetilde{\mathcal{N}}(x)$. We present
the absolute error value, as the approximation naturally oscillates around the true function. As $M$
increases, the point-to-point differences decrease between one and two orders of magnitude for every
additional ten mixands. Given enough components, we both approximate the distribution shape
accurately and achieve substantially smaller mixand standard deviations as compared to the original
Gaussian, by a factor of five in this example.

The lines in \cref{fig:GS_mism} show up to four possible distinct regions of approximation
behaviour: (i) the rapidly oscillating errors in the centre of the distribution, approximately for
$|x| \lesssim 2$ at large $M$, (ii) the comparatively smooth regions at the lower ends of the
Gaussian curve, $2 \lesssim |x| \lesssim 4 $, then (iii) somewhat increasing errors towards the
tails, and finally, (iv) no approximation beyond the outermost mixands (dashed lines with markers),
where the difference approaches the original function $\widetilde{\mathcal{N}}(x)$; the heuristic
approximation of \cref{sec:Heur:Odd,sec:Heur:Even} minimises the mismatch in this outer region. The
errors are larger near the centre of the plot, at the peak of the original Gaussian, and at the very
ends, where the outer-most mixands dissipate. The observed modes in the error behaviour might
indicate that although the proposed approximation is quite
accurate, it is still suboptimal. A non-uniform choice of the means and variances may provide
more uniform and lower errors, likely at the cost of increased complexity in finding such an
approximation.

\begin{figure}[htb!]
    \centering
    \includegraphics[width=\textwidth]{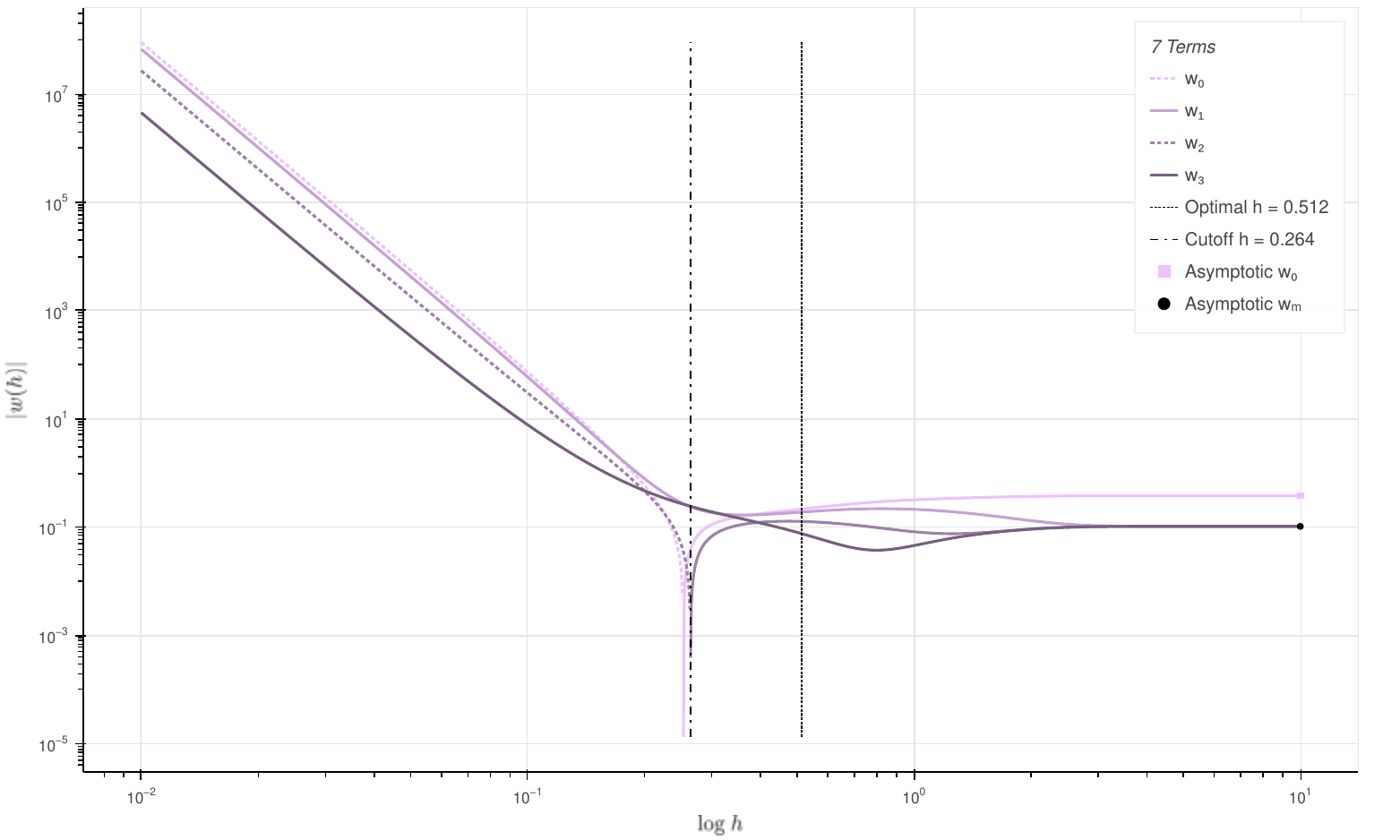}
    \caption{Dependence of the mixand weights on $h$ for $M = 7$. Dash-dotted vertical line:
    $h_\text{cut}$. Dotted vertical line: $h_\text{opt}$. Dot: asymptotic value for $w_m$, $m
    > 0$ as $h \to \infty$. Square: asymptotic value for $w_0$ as $h \to \infty$.}
    \label{fig:logw_M7}
    \vspace{-3mm}
\end{figure}

Moving next to the mixand weights, \cref{fig:logw_M7} displays the dependence of $|w_m|$ on the mean
step $h$ for $M = 7$; only the mixands with non-negative indices are used, as the rest have the same
weights by symmetry. The weights are plotted as solid lines when positive and as dashed lines when
negative. In this example, the weights $w_1$ and $w_3$ are always positive, while the weights $w_0$
and $w_2$ become negative for sufficiently small $h$. This pattern was observed consistently for all
considered numbers of mixands except $M = 2$, when the weights are always positive. All weights are
positive at $h_\text{opt}$, but then, as $h$ decreases, the weight of the second outer-most mixand
is the first to turn negative; for $M = 3$ this would be the weight $w_0$ of the central mixand. Heading
towards zero, the weights of the fourth, sixth, and so on, mixands as counted from the outer end
also turn negative, provided those mixands exist for the given $M$. Based on this pattern, we
compute the cutoff mean step $h_\text{cut}$ (dash-dot line in \cref{fig:logw_M7}) by numerically
finding the left-most root of the second outer-most weight. For $h < h_\text{cut}$ no solutions are
valid as at least some weights are negative. All weights
are positive to the right of $h_\text{cut}$ for small and intermediate values of $M$. However, the pattern of weight roots for very
large $M \gg 100$ is currently unknown.

On the left side of \cref{fig:logw_M7}, the weights asymptotically tend to infinity as $h \to 0$,
with ever-increasing rate as the number of mixands grows. For $M = 64$, the weights are by modulo on the
order of $10^{50}$--$10^{60}$ at $h = 0.01$. On the right end of the plot, we demonstrate that the
weights attain their expected asymptotic values for large $h$, in this case given by
\cref{eq:W:Odd:Inf}.

\begin{figure}[htb!]
    \centering
    \begin{subfigure}[htb]{\textwidth}
        \centering
        \includegraphics[width=\textwidth]{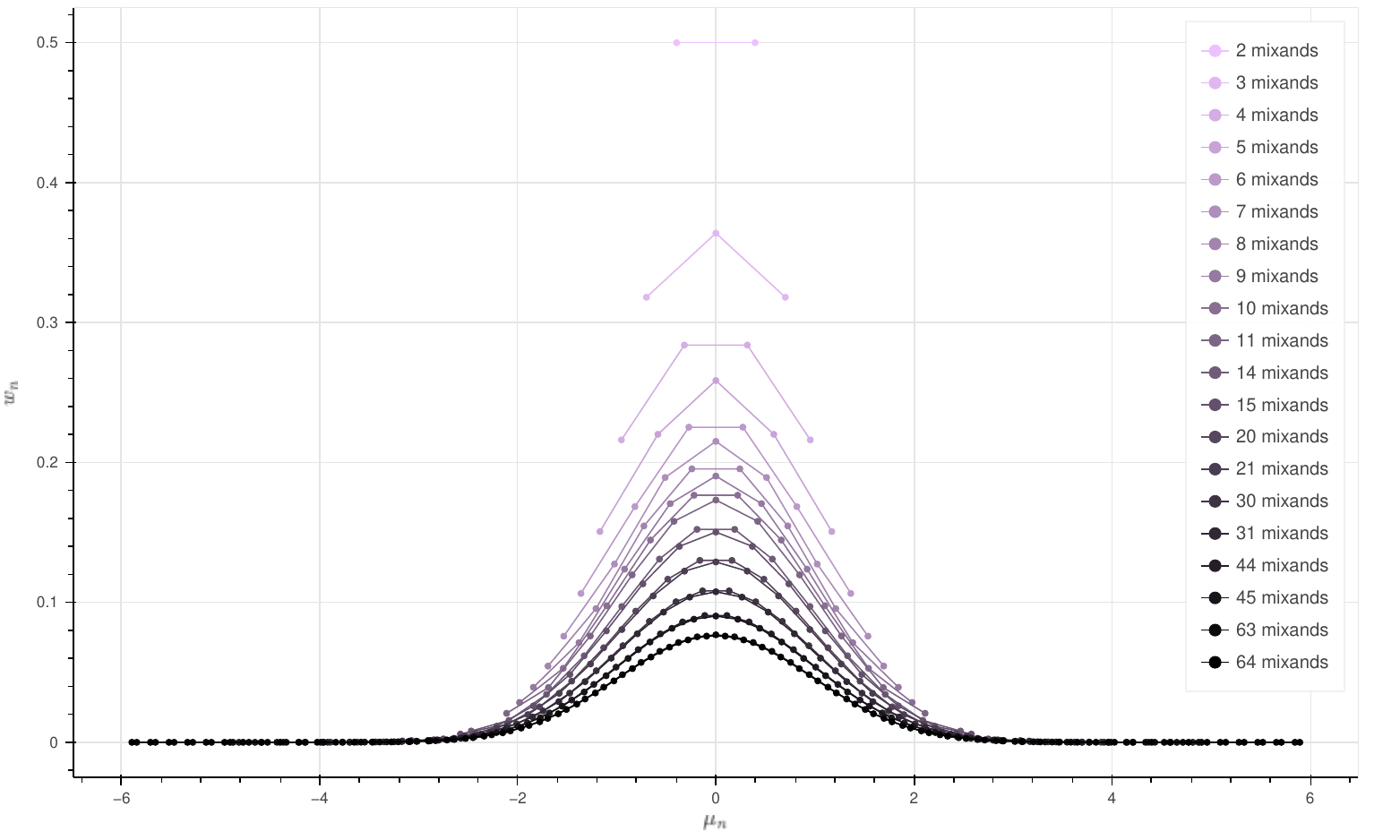}
        \caption{Dependence of the mixand weights on the position of their means.}
        \label{fig:wmu_M}
    \end{subfigure}
    \vfill
    \begin{subfigure}[b]{\textwidth}
        \centering
        \includegraphics[width=\textwidth]{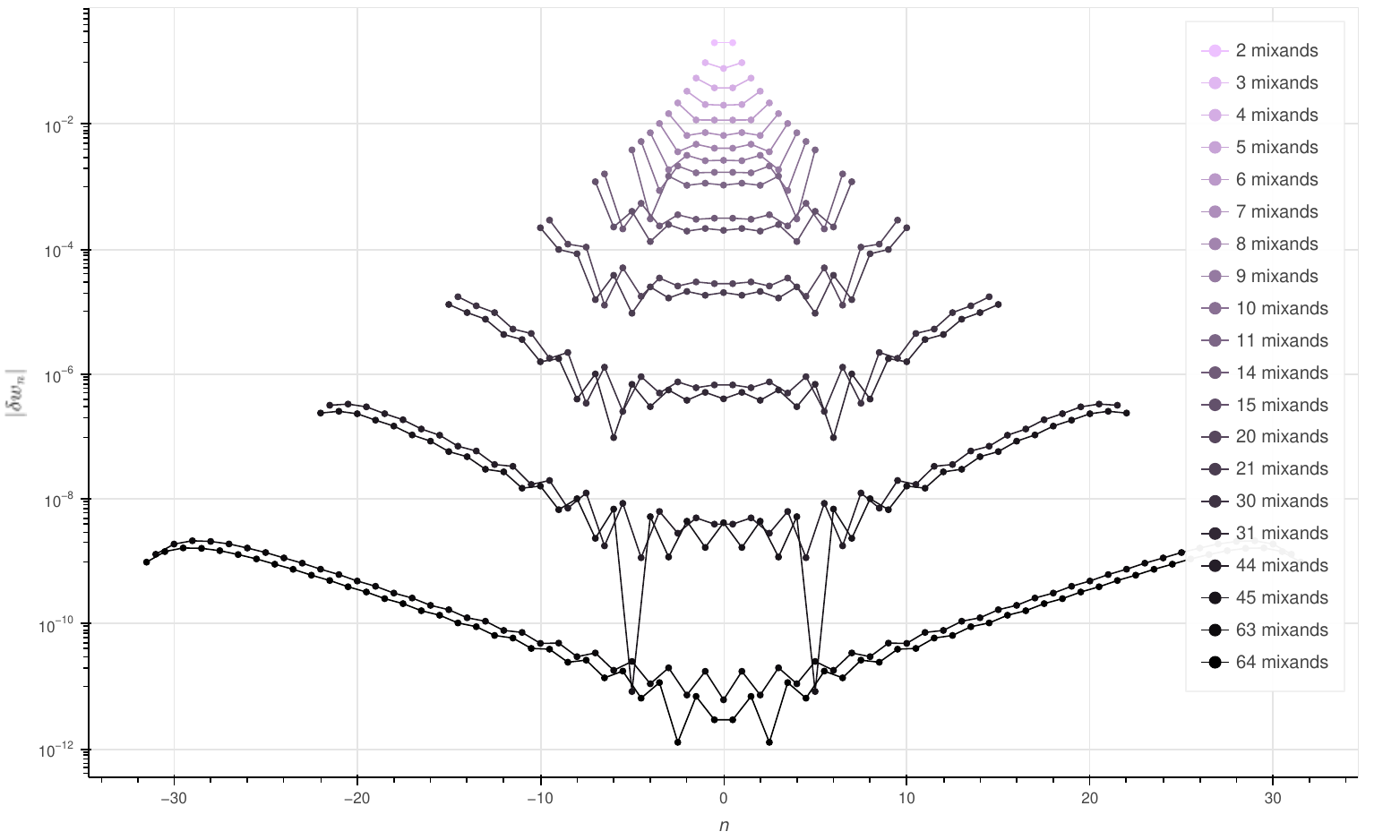}
        \caption{The absolute error in the approximate weights from \cref{eq:AsInfM} as compared to
        the exact solutions.}
        \label{fig:wmu_M_diff}
    \end{subfigure}
    \caption{The mixand weights (top) and the errors of the large-$M$ weight approximation
    (bottom).}
    \label{fig:weights}
    \vspace{-3mm}
\end{figure}

The overall layout of the mixand weights for various $M$ is presented in \cref{fig:wmu_M}. The
values on the $x$ axis are the corresponding positions of the mixand means. As demonstrated in
\cref{sec:LargeM:Odd,sec:LargeM:Even}, for large $M$, the weights' dependence on the off-centre
mixand index (and hence, the mean) approaches a Gaussian curve, which is confirmed by the
presented plot. \cref{fig:wmu_M_diff} displays the absolute difference between the precise weights
and their large-$M$ asymptotic approximations. The asymptotic solution is most accurate closer
to the centre of the original distribution, and, as expected, becomes more precise as $M$
increases.

\begin{figure}[htb!]
    \centering
    \includegraphics[width=\textwidth]{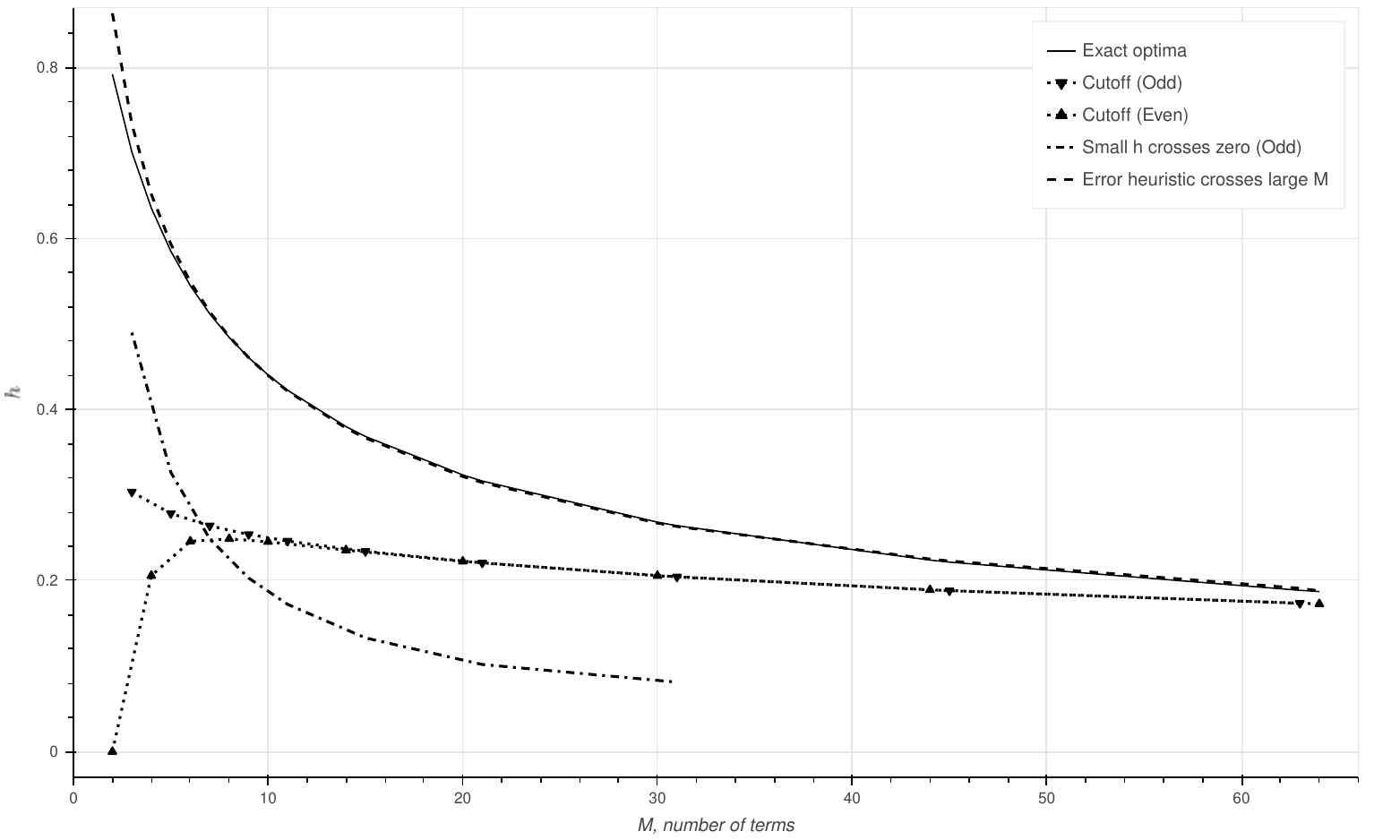}
    \caption{Relationship between $h$ and $M$ for the exact optima $h_\text{opt}$ (solid
    line), cutoff limits $h_\text{cut}$ (dotted lines with different triangle symbols for the
    odd and even cases), the intersection $h^\dagger_\text{opt}$ between the small-$h$
    approximation and zero (dash-dotted line, only up to $M = 31$ as the solution of
    \cref{eq:L2h0h2:Gsc} becomes numerically unstable beyond that value), and the intersection $h^*_\text{opt}$ between the
    heuristic solution and the large-$M$ asymptotic curve (dashed line).}
    \label{fig:hopt_M}
    \vspace{-3mm}
\end{figure}

\cref{fig:hopt_M} examines the dependence of the various special values of $h$
on the number of mixands $M$. The optimal mean step $h_\text{opt}$ (solid line)
decreases with $M$, as expected. However, it is unknown whether it converges to zero or to some
finite limit, which would be a function of the mixand standard deviation $\sigma$. This exact
solution is closely followed by the heuristic approximation $h^*_\text{opt}$ (dashed line);
however, being just a heuristic, it does not approach the true solution asymptotically. In fact, the
difference between the two exhibits a small increase closer to the right side of the plot. Still,
the agreement is good for the presented medium to large values of $M$. The small-$h$
approximation $h^\dagger_\text{opt}$, in comparison, is nowhere near the exact solution (dash-dots).
The dotted lines with triangle markers show the cutoff limits $h_\text{cut}$. Within the plot
limits, they are strictly below the optimal steps. However, it is not clear if these
lines would intersect for very large $M$.

\begin{figure}[htb!]
    \centering
    \includegraphics[width=\textwidth]{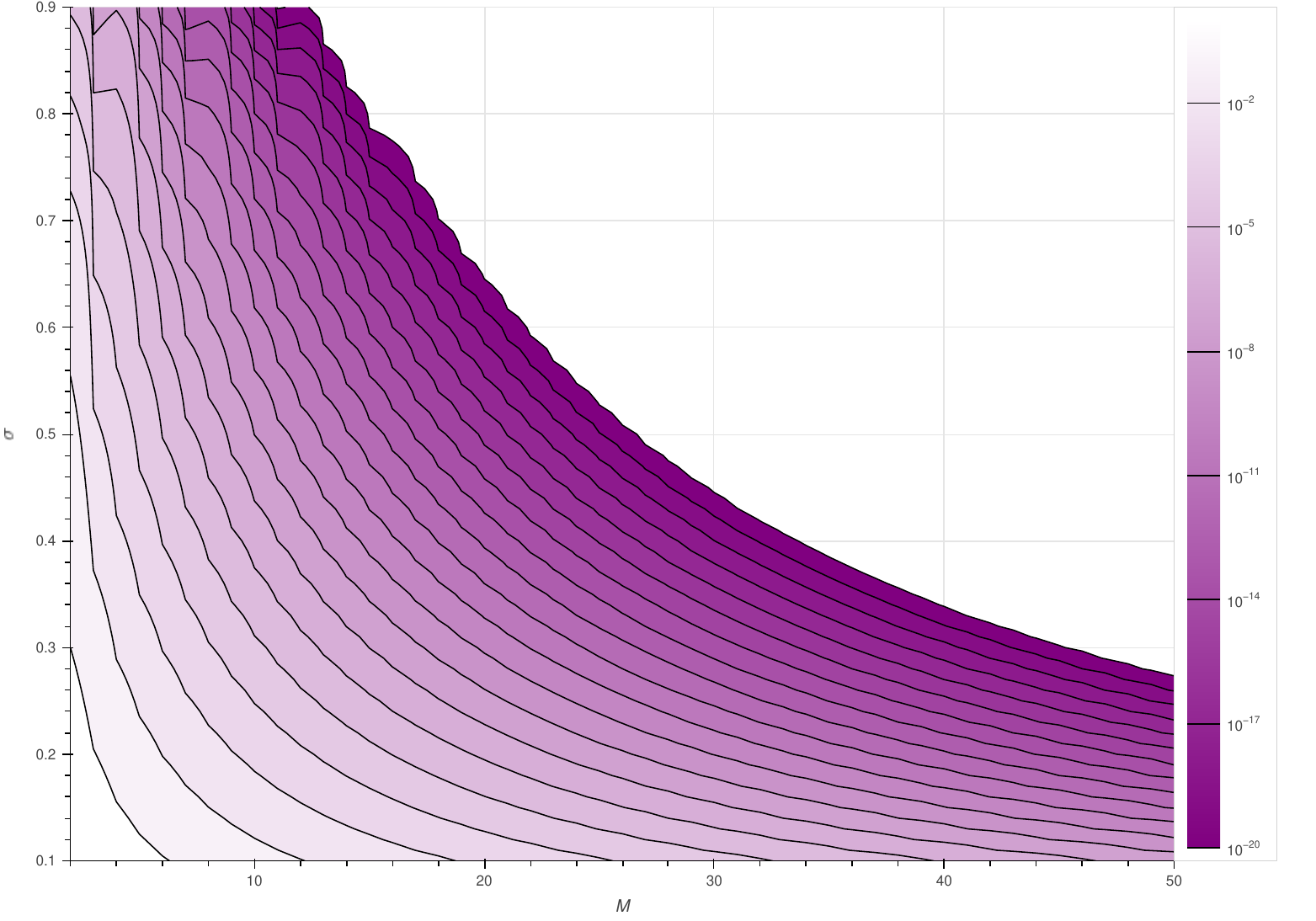}
    \caption{Relationship between $M$ ($x$-axis), $\sigma$ ($y$-axis) and $L^2$ at the point of
    optimum (colour-bar).}
    \label{fig:L2_sigma_M}
    \vspace{-3mm}
\end{figure}

So far, all examples assumed $\sigma = 0.2$ (although the presented conclusions
were confirmed for other values). In general, however, the choice of $\sigma$ is
made by the end-user based on the non-linearity of the problem at hand. Similarly, the user may have
some preferences for the number of terms $M$.
The heatmap in \cref{fig:L2_sigma_M}
shows the relationship between the optimal approximation accuracy $L^2$, the mixand variance $\sigma$, and $M$.
This plot enables the user to check the resulting
$L^2$ for the specified $\sigma$ and $M$, or, alternatively, select an appropriate $M$ based on a
chosen $\sigma$ and desired $L^2$. Naturally, smaller $\sigma$ values require larger $M$ to achieve
the same accuracy in $L^2$, whereas for $\sigma$ close to one, the original standard deviation, even
a few mixands provide highly accurate approximations. \cref{fig:L2_sigma_M} provides a
visual lookup-table for finding the right balance between these design parameters.

The top-right corner of \cref{fig:L2_sigma_M} ($\sigma$ close to one and large $M$), is filled with
white, because the $L^2$ values in this area are exceedingly small and fall out of the presented
colour spectrum limits. For reference, $L^2$ is about $10^{-60}$ in the top-right corner, and likely such an
extreme combination of $\sigma$, $M$, and accuracy has low practical importance. Thus, we
ignore such values instead of extending the colour spectrum to the full range of $L^2$, which
would result in lower resolution for moderate $L^2$, $M$, and $\sigma$, arguably of most interest in
applications.

\subsection{Evaluation}
\label{sec:Evaluation:1D}

As explained in \cref{sec:Recap}, our approach intends to preserve the distribution shape
while controlling mixand variance and quantity. Other considerations, for example, matching the
distribution moments, were so far explicitly ignored. However, we are approximating not just a
Gaussian function but a probability distribution, and therefore, should consider other measures
of approximation quality, arguably more statistical in nature than the $L^2$ mismatch.

\subsubsection{Moment comparison}
\label{sec:Moments:1D}

\begin{figure}[htb!]
    \centering
    \includegraphics[width=\textwidth]{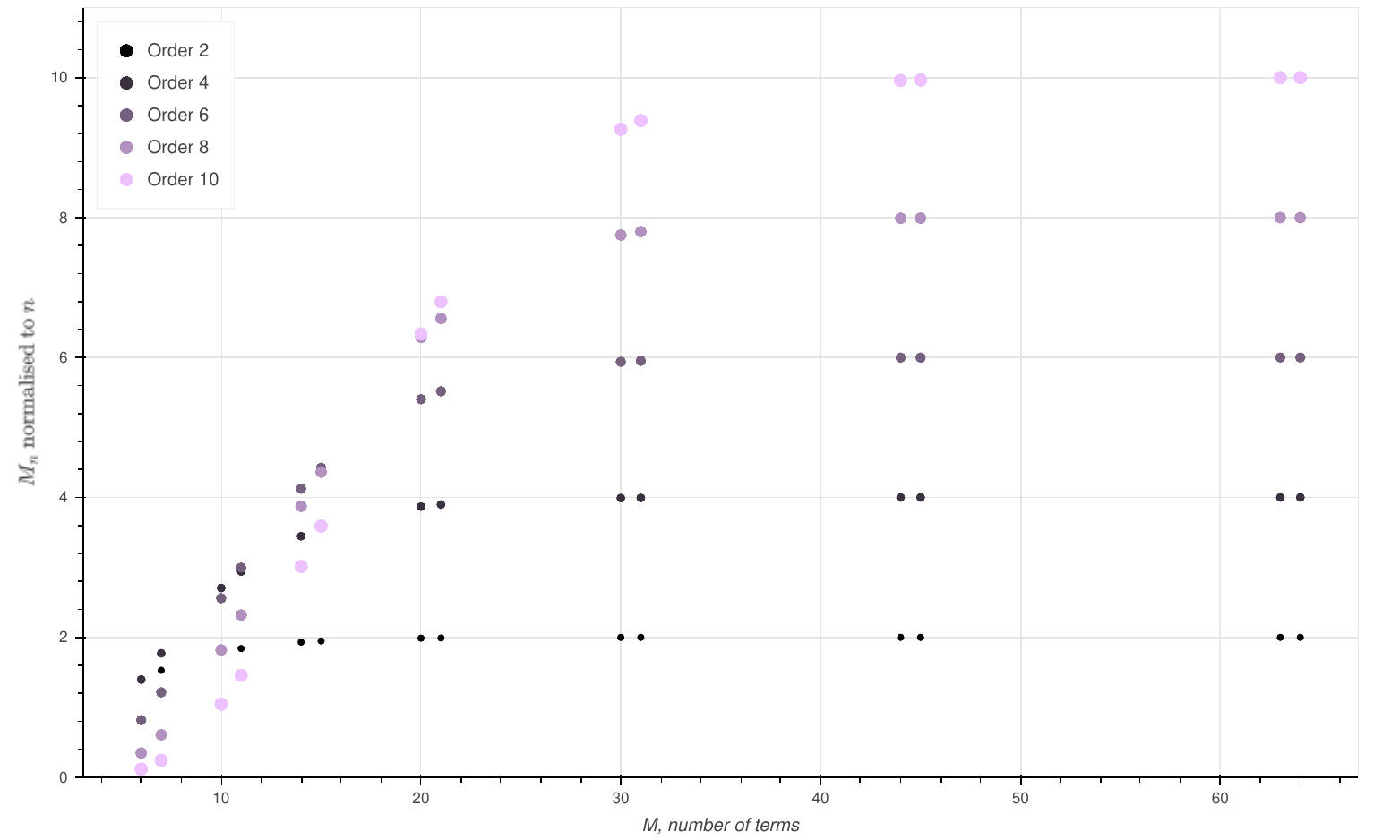}
    \caption{Moments of the Gaussian sum approximation normalised by \cref{eq:1d:Moment:Norm} for $n
    = 2$, $4$, $6$, $8$, and $10$. The shade and size of the scatter points indicates the moment
    order.}
    \label{fig:moms}
    \vspace{-3mm}
\end{figure}

\begin{figure}[htb!]
    \centering
    \includegraphics[width=\textwidth]{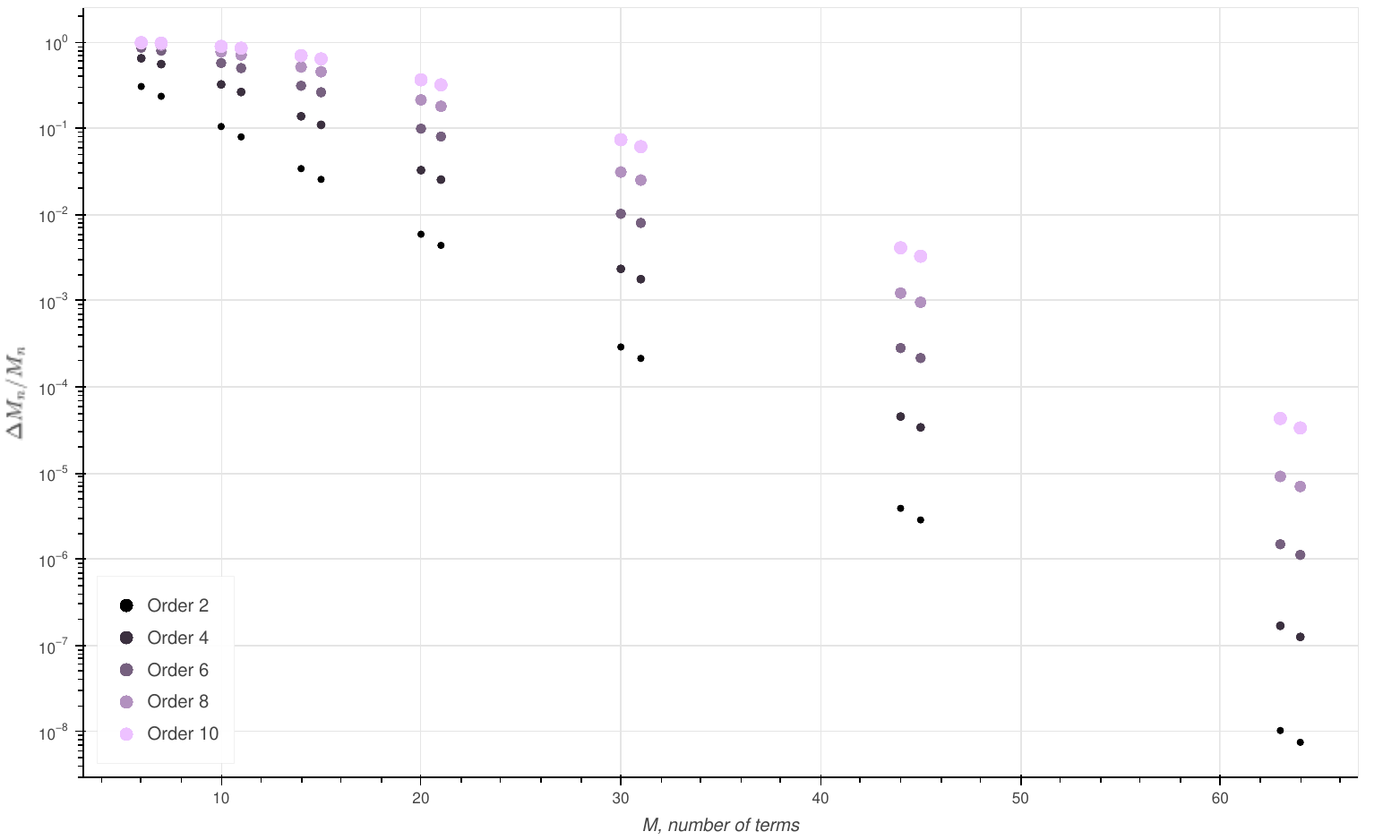}
    \caption{Relative error of moments, using the same shade and size conventions as in
    \cref{fig:moms}. Here $\Delta M_n = |\mathbb{E}[X^n] - M_n|$.}
    \label{fig:merrs}
    \vspace{-3mm}
\end{figure}

The proposed approximation does not use moment matching. Implicitly, we do match the odd-order
moments, as our solution is symmetric by design, and therefore, has all odd moments correctly equal
to zero. But none of the higher-order even moments are fixed. Still, we hope that matching the
distribution shape would result in \textit{some} approximation of those
moments. This supposition is examined in \cref{fig:moms,fig:merrs}.

Some special considerations are necessary to visualise moments of multiple orders at once, as
they rapidly grow with order \cite[p.~148]{Papoulis:2002:Probability}, and therefore, are hard to
see on the same scale. We normalise the moments for plotting as
\begin{equation}
    \label{eq:1d:Moment:Norm}
    \overline{M}_n = n \frac{\mathbb{E}[X^n]}{M_n} \text{,}
\end{equation}
where $\overline{M}_n$ is the value to plot, $M_n$ is the true $n^{\text{th}}$-order moment of the
original Gaussian, and $\mathbb{E}[X^n]$ is the corresponding moment of the mixture
distribution.\footnote{The absolute moments of mixands were computed following
\cite{Winkelbauer:2012:Moments}.}

\cref{fig:moms} demonstrates that as the number of mixands increases, the normalised moments
converge to their expected values, which are, by \cref{eq:1d:Moment:Norm}, the moment orders. The
five-fold decrease in the mixand standard deviation likely
precludes any good approximation for small $M$. Thus, we start with $M = 6$ and also skip some $M$
values, as compared to the previous figures, to de-clutter the plot. With $M$ growing, convergence
is achieved progressively for increasingly higher orders of moments: for the second moment between
$M = 15$ and $20$, for the fourth moment between $M = 20$ and $30$, and so on.

To quantify this initial visual assessment, \cref{fig:merrs} shows the relative errors of the
approximate moments with respect to the exact values. For a given $M$, the errors are consistently
higher for moments of higher order, as in \cref{fig:moms}. With increasing $M$, the
errors always decrease and, after some initial transition, become almost straight lines on the
logarithmic scale of the figure, \textit{i.e.} decrease almost exponentially. For the highest
order of approximation presented in the plot, $M = 64$, the distribution variance is approximated
with no more than $10^{-8}$ relative error and the $10^\text{th}$ order moment with less than
$10^{-4}$ error.

Thus, although our approximation technique never aimed at moment matching, it
actually performs this task well. It does not match \textit{any} moments precisely, but instead
approximates \textit{all} of them to some degree. An approach considering matching the moments
up to $10^\text{th}$ order would have to solve multivariate and likely non-polynomial equations
with variables to the power of at least the same $10^\text{th}$ degree, analyse multiple roots without
closed-form solutions, and possibly impose extra constraints to preserve positive mixand variances
(\textit{e.g.,} discussed by \cite{Leutnant:2011:Versatile}). This is decidedly more complicated
than solving linear systems of equations and univariate optimisation problems in
our method, and the numerical errors of such a solution may end up comparable in magnitude to our
accuracy.

\subsubsection{KL divergence}
\label{sec:KLD:1D}

As a second statistical measure of approximation quality we consider the Kullback-Leibler (KL)
divergence \cite{Cover:2006:Elements}
\begin{equation}
    \label{eq:KL:1D}
    D_{\text{KL}}(\widetilde{\mathcal{N}} || \widetilde{\mathcal{Q}})
        = \int_{-\infty}^{\infty} \widetilde{\mathcal{N}}(x)
                                      \log{\frac{\widetilde{\mathcal{N}}(x)}
                                                {\widetilde{\mathcal{Q}}(x)}} dx \text{.}
\end{equation}
The integral is evaluated numerically, which in our case is straightforward as the two
distributions are close and the appropriate finite integration limits follow from the explicit form
of $\widetilde{\mathcal{N}}$ and $\widetilde{\mathcal{Q}}$.

\begin{figure}[htb!]
    \centering
    \includegraphics[width=\textwidth]{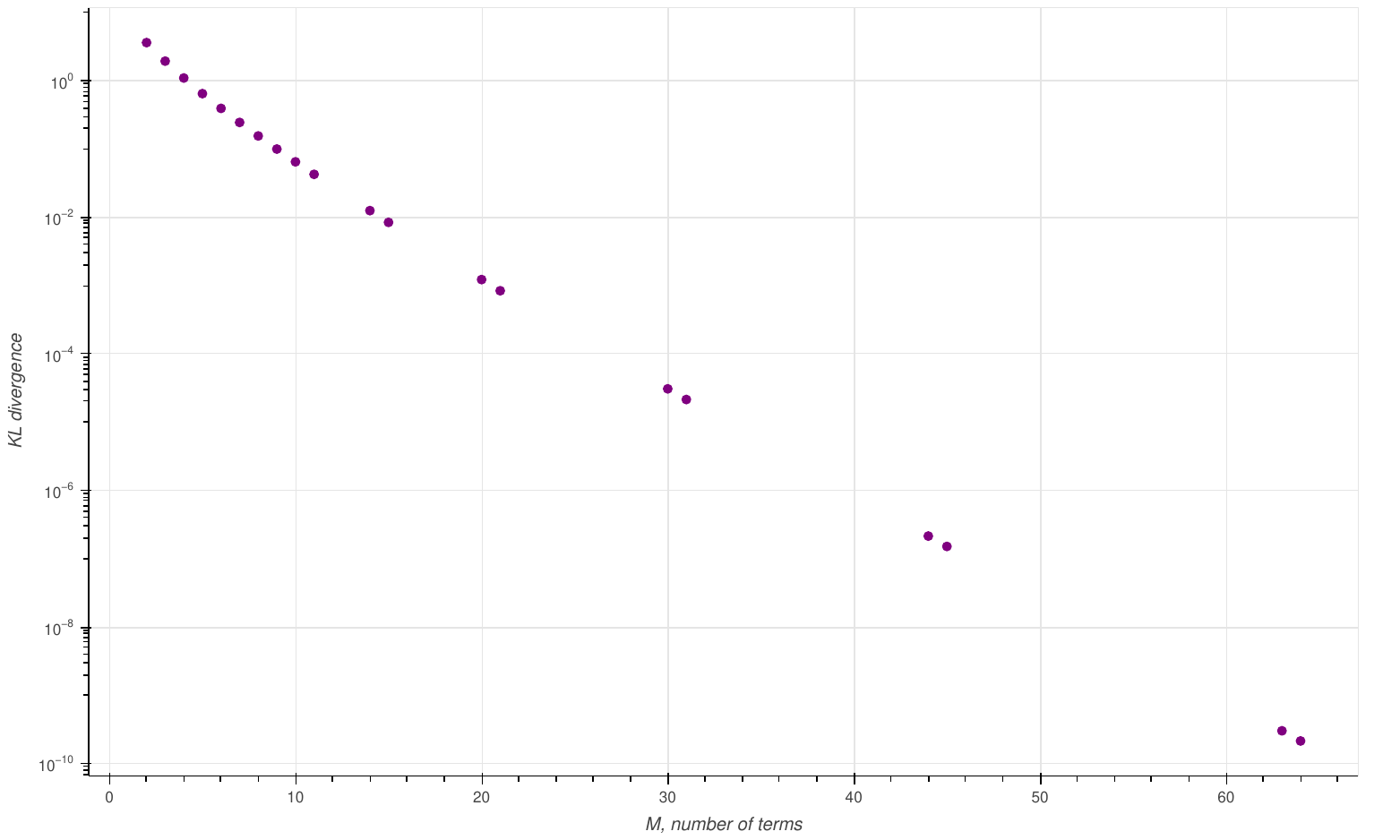}
    \caption{KL divergence between the original Gaussian and its approximations.}
    \label{fig:KL:1D}
    \vspace{-3mm}
\end{figure}

The resulting dependence of $D_{\text{KL}}(\widetilde{\mathcal{N}} || \widetilde{\mathcal{Q}})$ on
the number of mixands $M$ is shown in \cref{fig:KL:1D}. Similar to \cref{fig:merrs}, after some
initial transition range, the divergence exhibits nearly exponential decrease with $M$.

\section{Discussion\label{sec:Discussion}}

This paper presents an algorithm for approximating a Gaussian function with Gaussian Mixtures. The
obtained solution exhibits attractive properties. It is flexible, allowing the user to select the
right balance between the mixand variance, their number, not limited to low values, and the
approximation error (\cref{fig:L2_sigma_M}). The algorithm is robust: we demonstrated the
results for a wide range of input parameters. Most importantly, the approximation is highly
accurate and shows rapid convergence with increasing number of mixands for all the presented
measures of difference: the $L^2$ norm (\cref{fig:L2_sigma_M}), the point-to-point mismatch
(\cref{fig:GS_mism}), the moment errors (\cref{fig:merrs}) and the KL divergence (\cref{fig:KL:1D}).

Although the algorithm is sufficiently fast, it is best applied as a database or a
look-up library. The approximation is defined by two keys, usually the number of mixands $M$ and the desired
variance $\sigma^2$. Alternatively, one may specify the desired approximation accuracy $L^2_\text{max}$
and $M$ to find the minimal $\sigma$ that satisfies such requirements, or even use
$L^2_\text{max}$ and $\sigma$ to find the minimal $M$ (see \cref{fig:L2_sigma_M}). One of these keys
is already discrete and the other two can be sampled on a discrete grid; then all the approximation
parameters could be pre-computed, including the moment errors and KL divergences that can serve as
auxiliary search constraints. Storage requirements are modest, thus, grid sampling can be dense. The
splitting for a non-standard Gaussian is obtained by linear scaling and offset of the mixands.

As an alternative to the exact solution, for reasonably large $M$, such as $M \geq 15$, we can use the two closed-form solutions
obtained for special cases of input parameters: the heuristic method of
\crefrange{eq:Heur:Odd}{eq:Heur:Alt:Odd} and the large-$M$ asymptotic solution by \cref{eq:L2:Odd:AsympM}.
Their combination gives a good approximation to the exact minimum of $L^2$ (\cref{fig:L2h_M_hr}).
A simplified splitting algorithm can
\begin{itemize}
    \item Set the mean step $h = h^*_\text{opt}$, the intersection of the two special solutions, and
    \item Use this $h$ for the large-$M$ asymptotic weights from
        \cref{eq:AsInfM:Odd,eq:AsInfM:Even},
\end{itemize}
thus avoiding the need for the high-precision arithmetic, matrix solvers, and optimisers.
Still, the solutions by \cref{eq:Heur:Odd,eq:Heur:Alt:Odd} are only heuristic
approximations. Their accuracy has been established numerically for moderate-to-large values of $M$,
which are arguably most practically important, but the accuracy for
arbitrary $M$ is not guaranteed. In particular, this approach cannot be used to examine the asymptotic behaviour of
$h_\text{opt}$ for large $M$, which remains unknown.

The approximations with large $M$ open the possibility of accurate recursive splitting. If we know that
non-linearity is most pronounced near a particular point $x = x_\text{NL}$, we can split the
original Gaussian into a large but still moderate number of components, and then apply the same
algorithm to split one or several components nearest to $x_\text{NL}$.
For example, split into 10 mixands with moderate $\sigma$, thus achieving high approximation
accuracy, and then split one of those into 10 second-order components, ending up with 19 mixands in
total, 10 of which with a smaller $\sigma_\text{sub} = \sigma^2$. This is more efficient than
splitting into 100 $\sigma_\text{sub}$-wide components, as such resolution is
not needed everywhere.

The obtained GS is, of course, an approximation, as it is not possible to match a Gaussian exactly
using a limited number of mixands. However, by selecting the number of terms $M$ on demand, we
can drive the mismatch down to be within other sources of errors, \textit{e.g.}, the subsequent UKF
or CKF transformations in tracking applications, or even finite precision of floating-point numbers.
Our algorithm occupies the middle ground and provides an uninterrupted transition between low-term
techniques, as in \cite{Faubel:2010:Further, Hanebeck:2003:Progressive, DeMars:2013:Entropy},
and the generic large-$M$ Gaussian interpolation algorithms, such as \cite{Mazya:1996:Approximate}.
Compared to \cite{Mazya:1996:Approximate}, we use different starting requirements, but end up with
strikingly similar results.
The large-$M$ asymptotic of our solution
(\cref{sec:LargeM:Odd,sec:LargeM:Even}) reproduces the form of the quasi-interpolation formula from
\cite{Mazya:1996:Approximate}. Our solution is not valid for small steps $h$ between the component
means, while the quasi-interpolation from \cite{Mazya:1996:Approximate} does not converge for
small $h$. Convergence is demonstrated for some
schemes \cite{Beatson:1992:Quasi} with variable $\sigma$ decreasing with $h$, which for our case may
indicate that very large values of $M$ require correspondingly smaller $\sigma$ for the solution to
be valid.

The achieved high approximation accuracy over a range of design parameters and detailed partitioning
of the original distribution make our approximation a powerful tool. It can be applied to examine
the distribution at different scales, to subject it to highly non-linear transformations without a
loss of information, or to perform numerical integration. These and other applications are discussed
in \cite{Mazya:2007:Approximate} for their large-$M$ technique, and many of those equally apply to
our algorithm.

The efficiency and robustness of the presented method derive, in our opinion, from the selected
parametrisation of the mixture and the cost function. When searching for the
best approximation of a Gaussian with a GS, the optimisation variables in the order of
difficulty, from the easiest to the hardest, are arguably the weights, the means, and the variances.
We directly take the user input for the desired mixand variance, the same for all components,
so the most ``difficult'' variables in our ranking are fixed early. Then, the
position of the mixand means is described by a single parameter, the step $h$.
Finally, the cost function with quadratic dependency on
the weights leads to a system of linear equations for these variables. Together these design choices
reduce the original multivariate problem to univariate optimisation, irrespective of the number of
mixands. Even the univariate cost function is expected to have a reasonable behaviour on the
remaining unknown $h$. Too large steps leave gaps between the components, while steps that are too small
concentrate all components near zero and leave the tails poorly matched; both extremes should result
in larger mismatches. Intuitively, we expect no oscillatory dependence of $L^2$ on $h$, and
therefore, $L^2(h)$ to have a single minimum, as confirmed by
\crefrange{fig:L2h_M_c}{fig:L2h_M_hr}.

The approximation errors for the given $M$ can be further decreased and made more uniform (see
\cref{fig:GS_mism}) by extending the parametrisation approach, for example, using a two-parameter
description of the separation between the means.
Any two-argument (or even
$K$-argument) parameterisation of $h$ and / or $\sigma$ would still lead to a linear system of
equations for the weights. Solving that system would yield a bivariate (accordingly, $K$-variate)
form of $L^2$, which is easier to minimise than the GS mismatch in terms of the original variables.

Our parametric approach should also work for creating GS approximations of other continuous
distributions. A nearly closed-form solution is possible if the products $\langle \Phi,
\mathcal{N}_k(\sigma_k) \rangle$, where $\Phi(x)$ is the PDF function of that distribution, $k$ is the
component index and $\sigma_k^2$ its variance, are available in
closed form.
The optimal form of GS parametrisation may depend on the distribution.
In some cases, a practical solution may compute these integrals numerically.
If the
resulting approximation is applied as a library of pre-computed splittings, the integration
performance does not matter.
The Gaussian splitting proposed here can be used for the integration. Replacing
$\mathcal{N}_k$ by a GS of $\mathcal{N}_{k,m}$ with $\sigma_{k,m} \ll \sigma_k$, we
can evaluate each of $\langle \Phi, \mathcal{N}_{k,m} \rangle$ using a quadrature, which would work
the better the smaller $\sigma_{k,m}$ are.

One special case of ``solvable'' distributions is when $\Phi$ itself is a Gaussian Mixture. Then our
approach can be applied for Gaussian reduction. A similar technique was used in
\cite{Williams:2003:Cost}, but our method has a linear solution for the unknown weights.
However, the most important usage would be in multidimensional space,
where positions of the means are not obvious.

Most of the application examples that motivated this research are coming from the tracking
domain, which reflects our background. The algorithm itself is a general-purpose tool. Its
possible applications in image analysis \cite{Zhang:2003:EM}, general density estimation
\cite{Roeder:1997:Practical}, numerical analysis \cite{Mazya:2007:Approximate} or other areas are
beyond the intended scope of this study.

\section*{Acknowledgements}

The authors are deeply grateful to our colleagues Sanjeev Arulampalam, Kyle Talbot, Pei Leong, John
Maclean, Melissa Humphries, and Elena Kupriyanova for their valuable comments and suggestions to the
earlier versions of this manuscript. We highly appreciate the support of this work from our
employer, Acacia Systems, who are sponsoring the post-graduate studies of Athena and quietly
shouldered this major distraction by Dmitry.

\appendix
\section{Re-derivation of \cite{Leutnant:2011:Versatile} for $1$-D case\label{sec:Leutnant:1D}}

\renewcommand{\theequation}{A.\arabic{equation}}
\setcounter{equation}{0}

The splitting method of \cite{Leutnant:2011:Versatile} was obtained directly for the multivariate
problem.
Given our
focus on univariate distributions, we elucidate (and partially re-derive) their technique for
the univariate case.

In such restricted settings, there is only one possible ``direction of splitting'', and that is the
positive direction of the $x$ axis. In \cite{Leutnant:2011:Versatile} notation, this means setting
$\bm{u}_l \Rightarrow 1$ and $\widetilde{\bm{\Sigma}} \Rightarrow 1$ (the former direction vector
and the covariance matrix of the original distribution become scalars in one dimension). Then, from
Eq.~(6) of \cite{Leutnant:2011:Versatile} we obtain the variance of each mixand as
\begin{equation}
    \label{eq:LeutEq6}
    \sigma^2 = 1 - \beta \text{, with } \beta \defeq 2 \eta^2 \sum_{k = 1}^K w_k k^2 \text{,}
\end{equation}
where $\eta$ is the distance between the adjacent means of the components, assumed constant. For
brevity, we omit the lower index ``$0$'' for $\sigma$. The \textit{a priori} mixand weights from Eq.~(14) of
\cite{Leutnant:2011:Versatile} become
\begin{equation}
    \label{eq:LeutEq14}
    w_k = \frac{1}{C} \exp{\left( -\frac{1}{2} k^2 \eta^2 \right)} \text{,}
\end{equation}
where the normalisation ``constant'' is
\begin{equation}
    \label{eq:LeutEq14:C}
    C = \sum_{k = -K}^K \exp{\left( -\frac{1}{2} k^2 \eta^2 \right)} \text{.}
\end{equation}
We put ``constant'' in quotes because $C$, like $\beta$, is a function of $\eta$; for
simplicity, and following the discussion in \cite{Leutnant:2011:Versatile}, we set $\gamma = 1$ in
their Eq.~(14). Putting all these expressions together, we obtain
\begin{equation}
    \label{eq:LeutBeta}
    \beta = 2 \eta^2 \frac{ \sum_{k = 1}^K k^2 \exp{\left( -\frac{1}{2} k^2 \eta^2 \right)} }
                          { 1 + 2 \sum_{k = 1}^K \exp{\left( -\frac{1}{2} k^2 \eta^2 \right)} }
        \defeqr F(\eta^2, K) \text{.}
\end{equation}
Now, observe that $F(\eta^2, K)$ is a strictly increasing function of the number of approximation
terms $K$. Writing out only the numerator and omitting the $\eta^2$ argument of $F$ for brevity, we
get
\begin{equation*}
\begin{aligned}
    F(K + 1) - F(K) \!
        & \sim \! \Bigl( 1 + 2 \! \sum_{k = 1}^K \omega_k \Bigr) \! \sum_{k = 1}^{K + 1} k^2 \omega_k -
                  \Bigl( 1 + 2 \! \sum_{k = 1}^{K + 1} \omega_k \Bigr) \! \sum_{k = 1}^K k^2 \omega_k \\
        &= (1 + 2 S_0^K) \left( (K + 1)^2 \omega_{K + 1} + S_2^K \right) - \left(1 + 2S_0^K + 2\omega_{K + 1} \right) S_2^K \\
        &= (K + 1)^2 \omega_{K + 1} + 2 \omega_{K + 1} \left( (K + 1)^2 S_0^K - S_2^K \right) \\
        &= (K + 1)^2 \omega_{K + 1} + 2 \omega_{K + 1} \sum_{k = 1}^K \left( (K + 1)^2 - k^2 \right) \omega_k \text{,}
\end{aligned}
\end{equation*}
where we introduced $\omega_k = \exp{\left( -\frac{1}{2} k^2 \eta^2 \right)}$, $S_0^K = \sum_{k =
1}^K \omega_k$ and $S_2^K = \sum_{k = 1}^K k^2 \omega_k$. As $(K + 1)^2 - k^2 > 0$ for $k = 1,
\ldots, K$, the last line above is strictly positive, proving the statement. Therefore, the maximal
$\beta$ and the minimal component variance $\sigma^2$ are achieved when $K \to \infty$. Computing the exact
values for the infinite sum in \cref{eq:LeutBeta} is somewhat involved, so we use direct calculation for $K =
2^{16}$ to obtain $\beta \approx 0.926$ and $\sigma \approx 0.272$. Thus, for all the effort of creating
over 65,000 components, the standard deviation decreased by a factor of less than four.

\bibliographystyle{alpha}
\bibliography{gs-1D-full}

\end{document}